\documentclass[12pt, dvips, twoside]{amsart}
\usepackage{amsmath,amssymb,amscd,amsxtra,latexsym,%bbm,
graphicx,epsfig,%epic,eepic,oldgerm,bm,
euscript,%psfrag,
amsthm}
\usepackage{a4wide}
\usepackage[dvips]{color}
\theoremstyle{plain}
\newtheorem{theorem}{Theorem}[section]
\newtheorem*{theorem*}{Theorem}
\newtheorem{lemma}[theorem]{Lemma}
\newtheorem{proposition}[theorem]{Proposition}
\newtheorem{corollary}[theorem]{Corollary}

\theoremstyle{definition}  
\newtheorem*{remarks*}{Remarks}
\newtheorem{remark}[theorem]{Remark}

\newtheorem*{example*}{Example}
\newtheorem*{examples*}{Examples}
\newtheorem{definition}[theorem]{Definition}
\newtheorem*{definition*}{Definition}

\newcommand{\proofend}{\hspace*{\fill} $\Box$\\}

\def\1{\:\!}
\def\2{\;\!}

\def\Im{\operatorname {Im}}

\def\Diffc0{\operatorname{Diff^c_0}}

\def\Sympc0{\operatorname{Symp^c_0}}

\def\gl{\lambda}

\def\pp{\partial}

\def\ni{\noindent}
\def\b{\bigskip}

\def\.{\mskip1mu}
\def\?{\mskip-1mu}

\begin{document}

\title{Geometry and Topology of some overdetermined elliptic problems}

\author{Antonio Ros}
\address{(A.~Ros) 
Departamento de Geometr\'ia y Topolog\'ia,
Universidad de Granada,
Campus Fuentenueva,
18071 Granada,
Spain} 
\email{aros@ugr.es}

\author{Pieralberto Sicbaldi}
\address{(P.~Sicbaldi)
Laboratoire d'Analyse Topologie Probabilit\'es,
Universit\'e d'Aix-Marseille
39 Rue Joliot-Curie, 
13453 Marseille cedex 13,
France}
\email{pieralberto.sicbaldi@univ-cezanne.fr}

\date{\today}

\maketitle
\begin{abstract} 
We study necessary conditions on the geometry and the topology of domains in $\mathbb{R}^2$ that support a positive solution to a classical overdetermined elliptic problem
\[
\left\{\begin{array} {ll}
\Delta u + f(u) = 0 & \mbox{in }\; \Omega\\
               u= 0 & \mbox{on }\; \pp \Omega \\
\langle \nabla u, \nu \rangle = \textnormal{constant} &\mbox{on }\; \pp \Omega 
\end{array}\right.
\]
The ideas and tools we use come from constant mean curvature surface theory. In particular, we obtain a partial answer to a question posed by H. Berestycki, L. Caffarelli and L. Nirenberg in 1997. We investigate also some boundedness properties of the solution $u$. Some of our results generalize to higher dimensions. 
\end{abstract}

\section{Introduction}  \label{s:intro}

Let $\Omega$ be a connected open domain in $\mathbb{R}^n$ and $\overline{\Omega} =\Omega\cup\partial\Omega$. A long-standing open problem is to find necessary conditions on the geometry and the topology of $\Omega$ in order that the overdetermined elliptic problem 
\begin{equation}\label{pr_bis}
\left\{\begin{array} {ll}
\Delta u + f(u) = 0 & \mbox{in }\; \Omega\\
						u > 0 & \mbox{in }\; \Omega\\
               u= 0 & \mbox{on }\; \pp \Omega \\
\langle \nabla u, \nu \rangle = \alpha &\mbox{on }\; \pp \Omega 
\end{array}\right.
\end{equation}
is solvable, where $f$ is a given Lipschitz function, $\alpha$ is a nonpositive constant, $\nu$ is the unit normal vector about $\partial \Omega$, and $\langle \cdot, \cdot \rangle$ denotes the scalar product in $\mathbb{R}^n$.

\medskip

If $\Omega$ is bounded and regular enough, then the problem is very well understood: in a very interesting paper, \cite{Serrin}, J. Serrin proved that if $\Omega$ is a bounded domain, with boundary of class $C^2$, where there exists a solution $u \in C^2(\overline\Omega)$ to problem (\ref{pr_bis}), then $\Omega$ must be a ball. The proof of J. Serrin can be generalized in order to obtain the same result when $f$ is supposed to have only Lipschitz regularity, see \cite{Puc-Ser}. An alternative striking proof of some of the results of J. Serrin was given in \cite{Wei}. The result of J. Serrin has been of outstanding importance for two reasons: for applications to physics and to applied mathematics and for the development of very fruitful mathematical ideas. Indeed, problem (\ref{pr_bis}), when $f$ is constant, describes a viscous incompressible fluid moving in straight parallel streamlines through a straight pipe of given cross sectional form $\Omega$. If we fix rectangular coordinates $(x,y,z)$ with the $z$-axis directed along the pipe, it is well known that the flow velocity $u$ along the pipe is then a function of $x$ and $y$, and satisfies 
\[
\Delta u + k = 0
\]
where $k$ is a constant related to the viscosity and density of the fluid. The adherence condition is given by $u=0$ on $\partial \Omega$. The result of J. Serrin allows us to state that the tangential stress per unit area on the pipe wall (represented by $\mu\, \langle \nabla u, \nu \rangle$, where $\mu$ is the viscosity) is the same at all points of the wall if and only if the pipe has a circular cross section.
Other models from physics are also referable to problem (\ref{pr_bis}), for example in the linear theory of torsion of a solid straight bar of cross section $\Omega$, see \cite{Sok}. In this framework, the result of J. Serrin states that when a solid straight bar is subject to torsion, the magnitude of the resulting traction which occurs at the surface of the bar is independent of the position if and only if the bar has a circular cross section. Problem  (\ref{pr_bis}) is also related to a lower-dimensional obstacle problem (the so called Signorini problem, see \cite{Fre}). But besides the many applications, the paper by J. Serrin was very important because it made the moving plane method available to a large part of the mathematical community. This method had been introduced some years before by A. D. Alexandrov to prove that the only compact, constant mean curvature hypersurfaces embedded in $\mathbb{R}^n$ are the spheres, see~\cite{Alex}. The use of the moving plane method in analysis originated many fundamental results, such as the ones in \cite{GNN}.

\medskip

Overdetermined boundary conditions arise naturally also in free boundary problems, when the variational structure imposes suitable conditions on the separation interface: see for example \cite{Alt-Caf}. In this context it is important to underline that several methods for studying locally the regularity of solutions of free boundary problems are often based on blow-up techniques applied to the intersection of $\Omega$ with a small ball centered in a point of $\partial \Omega$, which lead then to the study of an elliptic problem in an unbounded domain. In this framework, problem (\ref{pr_bis}) in unbounded domains was considered by H. Berestycki, L. Caffarelli and L. Nirenberg in \cite{BCN}. In this case, the ÒtypicalÓ nonlinearity $f$ taken into account was $f(u) = u-u^3$, which reduces the equation in (\ref{pr_bis}) to the Allen-Cahn equation. Under the assumptions that $\Omega$ is a Lipschitz epigraph with some suitable control at infinity for its boundary, they proved that if problem (\ref{pr_bis}) admits a smooth, bounded solution, then $\Omega$ is a half-space. In the same paper, H. Berestycki, L. Caffarelli and L. Nirenberg proposed a very nice conjecture, that could be considered as the parallel of the result of J. Serrin that could be expected for overdetermined problems in domains not supposed to be bounded. The conjecture can be state as following:

\medskip

{\it Conjecture (C).} If $f$ is a Lipschitz function on $\mathbb{R}_+$ and $\Omega$ a smooth domain in~$\mathbb{R}^n$ such that $\mathbb{R}^n \backslash \overline{\Omega}$ is connected, then the existence of a bounded
solution to problem \eqref{pr_bis} implies that $\Omega$ is either a ball, a half-space, a generalized cylinder $B^k \times \mathbb{R}^{n-k}$ where $B^k$ is a ball in $\mathbb{R}^k$, or the complement of one of them.

\medskip

Such conjecture was motivated by the results obtained by H. Berestycki, L. Caffarelli and L. Nirenberg in \cite{BCN} and also by some results of W. Reichel in \cite{reichel} and A. Aftalion and J. Busca in \cite{aft-bus}. In these two last papers, authors were interested by overdetermined elliptic problems in  exterior domains, i.e. domains that are the complement of a compact region. Assuming that $\Omega$ is the exterior of some bounded and smooth region $D$ and $u$ is a bounded solution of problem (\ref{pr_bis}) for some particular classes of function $f$ and with some assumptions on the behavior of $u$ at infinity, they proved that $D$ is a ball.

\medskip

In some recent papers, A. Farina and E. Valdinoci obtained natural assumptions under which one can conclude that if $\Omega$ is an epigraph where there exists a solution to problem (\ref{pr_bis}) then $\Omega$ must be a half-space and $u$ is a function of only one variable, see~\cite{Far-Val-0},~\cite{Far-Val-2} and~\cite{Far-Val-1}. 
%In these papers many results are proved, and we want to recall some of them. A first result they obtain is that if $\Omega$ is a globally Lipschitz smooth epigraph of $\mathbb{R}^n$, with $n \geq 2$, then there exists no solution $u \in C^2(\overline{\Omega}) \cap L^\infty(\Omega)$ for the problem \eqref{pr_bis} without the Neumann condition at the boundary. In fact, the proof of Farina and Valdinoci applies to prove the same result for problem  \eqref{pr_bis} when $f(u) \geq \lambda\,u$ for a given positive constant $\lambda$. A second result of A. Farina and E. Valdinoci, that we will use in our paper, is the following (we refer to it as Theorem 1.6 of \cite{Far-Val-0}, even though this is only a part of the original theorem): if $f$ is locally Lipschitz, $n=2$, and $\Omega$ is a uniformly Lipschitz coercive epigraph of class $C^3$, then there exists no function $u \in C^2(\overline{\Omega}) \cap L^\infty(\Omega)$ which is solution of (\ref{pr_bis}) (the definition of uniformly Lipschitz coercive epigraph is presented below). An other result of A. Farina and E. Valdinoci, that also we will use in our paper, is the following (we refer to it as Theorem 1.2 of \cite{Far-Val-0}, even though this is only a part of the original theorem): if $f$ is locally Lipschitz, $n=2$, and $u \in C^2(\overline{\Omega})$ is a solution of (\ref{pr_bis}) increasing in the second variable and with bounded gradient, then $\Omega$ is a half-plane.

\medskip

In \cite{Sicbaldi}, P. Sicbaldi provided a counterexample to Conjecture (C) in dimension bigger or equal then 3 when $f$ is the linear function $f(t) = \lambda\, t$, $\lambda > 0$, constructing a periodic perturbation of the straight cylinder $B_1^n \times \mathbb{R}$, where $B_1^n$ is the unit ball of $\mathbb{R}^n$, that supports a periodic (and then bounded) solution to problem (\ref{pr_bis}).  In \cite{Sch-Sic}, F. Schlenk and P. Sicbaldi improved such result by constructing a smooth 1-parameter family of unbounded domains $s \mapsto \Omega_s$ in $\mathbb{R}^{n+1}$ for $n \geq 1$, whose boundaries are smooth periodic hypersurfaces of revolution with respect to an $\mathbb{R}$-axis 
and such that~\eqref{pr_bis} has a bounded solution in $\Omega_s$. They proves this result by showing that the cylinder $B_1^{n} \times \mathbb{R} \subset \mathbb{R}^{n+1}$ 
(for which it is easy to find a bounded solution to~\eqref{pr_bis} with $f(t)= \lambda\,t$)
bifurcates into unbounded domains whose boundary is a periodic 
hypersurface of revolution with respect to the axis of the cylinder,
and such that (\ref{pr_bis}) has a bounded solution with $f(t) = \lambda\, t$, $\lambda > 0$. The technique used by F. Schlenk and P. Sicbaldi is based on the Crandall-Rabinowitz bifurcation theorem, see \cite{cran-rab}, and the very interesting aspect of this result is that it is obtained by paralleling in a very strong sense the construction of the well known family of constant mean curvature surfaces of Delaunay in $\mathbb{R}^3$, which can easily be obtained by bifurcation from a straight cylinder with a bifurcation result as the Crandall-Rabinowitz theorem. The existence of the domains $\Omega_s$ provides a smooth 1-parameter family of counterexamples to Conjecture (C) in dimension bigger or equal then 3, but not in dimension 2 because in this case, $\Omega_s$ is a perturbation of a strip in $\mathbb{R}^{2}$ and then its complement is not connected. In dimension 2 Conjecture (C) is still open and it will be the starting point of this paper.

\medskip

Problem (\ref{pr_bis}) in the interesting case when $f = 0$ has been recently studied in $\mathbb{R}^2$ by F. H\'el\`ein, L. Hauswirth and F. Pacard. In \cite{HHP} they provide the following nontrivial example of domain where problem (\ref{pr_bis}) can be solved with $f = 0$:
\[
\Omega_{*} = \left\{\omega \in \mathbb{C}\,\,:\,\, |\Im \omega| < \frac{\pi}{2} + \cosh(\Re \omega) \right\},
\]
and they conjecture that $\Omega_{*}$, the half-planes and the complements of a ball are the only domains in $\mathbb{R}^2$ where problem (\ref{pr_bis}) with $f = 0$ can be solved. It is important to underline that their work is inspired by the theory of minimal surfaces, and it is interesting to remark that domains where problem (\ref{pr_bis}) with $f = 0$ can be solved arise as limits under scaling of sequences of domains where problem (\ref{pr_bis}) with $f(t) = \lambda\, t$ can be solved, just like minimal surfaces arise as limits under scaling of sequences of constant mean 
curvature surfaces.

\medskip

The result of J. Serrin in \cite{Serrin}, the result of P. Sicbaldi and F. Schlenk in \cite{Sch-Sic} and the result of F. H\'el\`ein, L. Hauswirth and F. Pacard in \cite{HHP} show that the geometry of overdetermined elliptic problems shares profound similarities with the theory of constant mean curvature hypersurfaces, even though the link between the two objects is not clear. Based on this facts, in this paper we want to study overdetermined problem using ideas and tools coming from the theory of constant mean curvature hypersurfaces. This strategy was not really exploited in the past, when in general overdetermined problems were studied with tools coming from PDEs and Analysis theories. One of our main results is the following:

\begin{theorem} \label{Tmain} In dimension $n=2$, the conjecture of Berestycki-Caffarelli-Nirenberg (C) is true in the following two cases:
\begin{itemize}
\item[(A)] when $\Omega$ is contained in a half-plane and $|\nabla u|$ is bounded, or
\item[(B)] when there exists a positive constant $\lambda$ such that $f(t) \geq \lambda\, t$, $\forall t>0$.
\end{itemize}
\end{theorem}

Note that in the case (A) we do not assume any hypothesis on the function $f$ and in (B) we do not assume any condition on the domain $\Omega$. Moreover in the case (A) the hypothesis that $|\nabla u|$ is bounded can be removed in some geometric situations (see Remark \ref{nablau_acotado} and Statement (T8) of Theorem \ref{T4}).

\medskip

In fact, in this paper we are interested in understanding the geometry and the topology of Euclidean 2-dimensional domains where problem (\ref{pr_bis}) can be solved, starting from the persuasion that such domains shares many properties with 2-dimensional constant mean curvatures surfaces in $\mathbb{R}^3$. In this framework, we will start from some very known and classical properties of constant mean curvatures surfaces in $\mathbb{R}^3$ and we will try to adapt its to overdetermined elliptic problems. Theorem \ref{Tmain} will be a corollary of some more general results, that we will present in the next section. 

\medskip

As a final remark, it is important to underline that the study of the geometry and topology of the elliptic overdetermined problems (\ref{pr_bis}) can be useful in order to develop new ideas which can be exploited in the understanding of the Schiffer's conjecture, that can be stated as following: if $\Omega$ is a bounded domain in $\mathbb{R}^n$ and there exists a solution $u \in C^2(\overline\Omega)$ of
\begin{equation}\label{shiffer}
\left\{\begin{array} {ll}
\Delta u + \lambda\, u = 0 & \mbox{in }\; \Omega\\
               u= \textnormal{constant} & \mbox{on }\; \pp \Omega \\
\langle \nabla u, \nu \rangle = \textnormal{0} &\mbox{on }\; \pp \Omega 
\end{array}\right.
\end{equation}
for some constant $\lambda$, then $\Omega$ is a ball (see \cite{Yau}). The study of the Schiffer's conjecture is now considered one of the outstanding problems in analysis since S. Williams proved in 1976 that the conjecture is equivalent to the famous Pompeiu problem in integral geometry, see \cite{Wil}. For a survey on this subject we remind to \cite{Zal}.

\b
\ni
{\bf Acknowledgments.}
This paper was started in April 2011, when the second
author visited the University of Granada. The first author is partially supported by MEC-FEDER Grants MTM2007-61775, MTM2011-22547 and J. Andaluc\'ia
Grant no. P09-FQM-5088. The second author is grateful to the members of
the Department of Geometry and Topology of the University of Granada for
their warm hospitality, to the International Scientific Coordination
Network `Geometric Analysis' (France and Spain) and to the applied analysis group of the LATP of Marseille for financial supporting. Both authors thank Alberto Farina and Joaqu\'in P\'erez for useful suggestions. 

\section{From constant mean curvature surface to overdetermined problems. Statement of the results}

The main aim of this paper is to study necessary geometric and topological conditions of domains $\Omega$ where the overdetermined elliptic problem (\ref{pr_bis}) can be solved, and we will focus our attention on 2-dimensional domains. It is convenient to give a name to such domains, and in this paper we will refer to them as {\it $f$-extremal domains}. The motivation of such definition comes from extremal domains for the first eigenvalue of the Laplacian. Indeed, 
if we consider the Dirichlet problem in a domain $\Omega$ of a Riemannian manifold $(M,g)$
\begin{equation} \label{e:sys1}
\left\{\begin{array} {ll}
\Delta_g\, u + \gl\, u = 0 &\mbox{in }\; \Omega\\
u=0 & \mbox{on }\; \pp \Omega 
\end{array}\right.
\end{equation}
where $\Delta_g$ is the Laplace-Beltrami operator (the natural generalization of the Euclidean Laplacian), we can denote by $\gl_1(\Omega)$ the smallest positive constant $\gl$ for which this system has a nonzero solution (i.e.,~$\gl_1(\Omega)$ is the first eigenvalue of the Laplace-Beltrami operator on $\Omega$ with 0 Dirichlet boundary condition). The solution $u$, up to a constant factor, is the only eigenfunction with constant sign in~$\Omega$ and we can consider~$u$ 
to be positive on~$\Omega$, see \cite{chavel}. Consider the functional $\Omega \to \gl_1 (\Omega)$ for all 
smooth bounded domains $\Omega$ in $M$ of the same volume, 
say $\textnormal{Vol}(\Omega) = V$. 
A classical result due to P. R. Garabedian and M. Schiffer in the Euclidean case, generalized by A. El Soufi and S. Ilias in the Riemannian one, asserts 
that $\Omega$ is a critical point for $\gl_1$ 
(among all domains of volume $V$) if and only if 
the first eigenfunction of the Laplace-Beltrami operator in $\Omega$ 
with $0$~Dirichlet boundary condition has also constant Neumann data 
at the boundary, see~\cite{Gar-Schif} and \cite{ElS-Il}.
In this case, $\Omega$ is called {\it extremal domain 
for the first eigenvalue of the Laplace-Beltrami operator} or simply {\it extremal domain}.
Extremal domains for the first eigenvalue of the Laplace-Beltrami operator are then characterized as the domains for which there exists a positive constant $\lambda$ such that the overdetermined elliptic problem
\begin{equation}\label{pr_1}
\left\{
\begin{array} {ll}
\Delta_g\, u + \gl\, u = 0 &\mbox{in }\; \Omega\\
u > 0  & \mbox{in }\; \Omega \\
u=0 & \mbox{on }\; \pp \Omega \\
g(\nabla u, \nu ) = \textnormal{constant} &\mbox{on }\; \pp \Omega
\end{array} 
\right.
\end{equation}
can be solved, where $\nu$ is the unit normal vector to $\pp \Omega$ pointing outwards $\Omega$.
If $\Omega$ is an unbounded domain of $\mathbb{R}^n$ the geometric meaning of extremal domain fails in general, except for the case when along each coordinate direction of $\mathbb{R}^n$ the domain $\Omega$ is bounded or periodic. In the case of periodic directions, one obtains extremal domains for the first eigenvalue of the Laplace-Beltrami operator in flat tori. An example of such domains are domains $\Omega_s$ found by F. Schlenk and P. Sicbaldi in \cite{Sch-Sic} and described in the previous section. Anyway, it is quite natural to continue to call a general domain (bounded or unbounded) of $\mathbb{R}^n$ where problem (\ref{pr_1}) can be solved an {\it extremal domain}. When we consider problem (\ref{pr_bis}) instead of  (\ref{pr_1}), we suggest to talk about {\it $f$-extremal domains}.

\medskip

We remark that if \eqref{pr_bis} is solvable, then the solution $u$ is unique, also for unbounded domains. In fact, if $u_1$ and $u_2$ are two functions that satisfy the elliptic equation of system \eqref{pr_bis} such that there exists a $C^1$ subset $\Gamma$ of $\partial \Omega$ where $u_1 = u_2$ and $\langle \nabla u_1, \nu \rangle = \langle \nabla u_2, \nu \rangle$, then $u_1=u_2$ in the whole $\Omega$, see~\cite{Far-Val-1}. It is clear that in \eqref{pr_bis}, the constant $\alpha$ must be non-positive because $u$ is positive in $\Omega$.

\medskip

As we said in the previous section, our work starts from the persuasion that constant mean curvature surfaces and $f$-extremal domains share profound similarities. For constant mean curvature surfaces in $\mathbb{R}^3$ the theory is very rich and many results are known. Our aim is to analyze the behavior of constant mean curvature surfaces in order to have new ideas on the behavior of extremal domains, and use these ideas in order to prove non-trivial results about overdetermined elliptic problems.  In this paper we are interested in three important results about constant mean curvature (hyper-)surfaces and for each of them we will show that there exists a parallel results about overdetermined elliptic problem. In order to state the first one, let us recall some topological facts. Denote by $B^{n}_{R} \subset \mathbb{R}^n$ the open ball centered at the origin with radius $R$. A properly embedded hypersurface $M$ in $\mathbb{R}^n$ is said to have {\it proper finite topology} if for a large $R$ we have that 
$M \backslash B^{n}_R$ has a finite number of connected noncompact components, the sphere $\partial B^n_{R}$ intersects $M$ transversally, and each component $E$ of $M \backslash \overline{B^{n}_R}$ is diffeomorphic to  $S^{n-1}_{1}\times [0,\infty[$, where $S^{n-1}_1$ is the boundary of $B^n_1$. Such component $E$ is called an {\it annular end} of $\Omega$. The first result we are interested in is the following very classical:

\begin{theorem} \label{Meeks_T} 
(W. H. Meeks, \cite{Meeks}). If $E$ is an annular end of a properly embedded, non-zero constant mean curvature surface $M$ in $\mathbb{R}^3$ of proper finite topology, then
\begin{enumerate}
	\item[(1)] $E$ stays at bounded distance from a straight line;
	\item[(2)] $M$ cannot have only one annular end;
	\item[(3)] If $M$ has exactly two annular ends, then $M$ stays at bounded distance from a straight line.
\end{enumerate}
\end{theorem}

The definition of annular end can be generalized in some sense for unbounded domains of $\mathbb{R}^n$. We can say that the domain $\Omega \subset \mathbb{R}^n$ has {\it
finite topology} if outside of a ball $B^n_R$ of large radius $R$, we
have that either
\begin{itemize}
\item $\overline{\Omega} \backslash \overline{B^{n}_R}$ is empty and then
$\overline{\Omega}$ is compact, or 
\item $\overline{\Omega} \backslash \overline{B^{n}_R}$ is equal to $\mathbb{R}^n \backslash
\overline{B^{n}_R}$ and then $\Omega$ is the complement of a compact region,
or 
\item $\overline{\Omega} \backslash B^{n}_R$ has a finite number of connected noncompact components and each
component $E$ is diffeomorphic to  $\overline{B^{n-1}_{1}}\times
[0,+\infty[$.
\end{itemize}
In the last case we can assume that the sphere $\partial
B^n_{R}$ intersects $\partial \Omega$ transversally and that each component of
$\partial B^n_{R}\cap\partial \Omega$ is diffeomorphic to
$\partial{B^{n-1}_{1}}$. Then, in the last case we will say that $\Omega$ has
{\it proper finite topology} and $E$ is a {\it solid
cylindrical end} of $\Omega$ if $n \geq 3$ or a {\it planar strip end} of $\Omega$ if $n=2$. If $n=2$, then $\Omega$ has proper finite
topology if and only if it is noncompact, $\partial \Omega$ has a finite
number of boundary components, some of them being
noncompact. Moreover the ends of $\Omega$
have the topology of a half-strip $[0,1]\times [0,+\infty[$. In this case,
the number of ends coincides with the number of noncompact components of
$\partial \Omega$.

\medskip

Let now define the following property:
\begin{enumerate}
\item [$P_1$] : there exists a positive constant $R$ such that $\overline{\Omega}$ does not contain any closed ball of radius $R$.
\end{enumerate}
%It is clear that property $P_{0}$ implies property $P_{1}$. 
Inspired by the result of W. H. Meeks, we will prove the following: 

\begin{theorem} \label{T2} Let $\Omega$ be an $f$-extremal domain of $\mathbb{R}^2$ of finite topology, satisfying the property $P_1$. Then, the following properties hold:
\begin{enumerate}
	\item[(T1)] If $E$ is a planar strip end of $\Omega$, then $E$ stays at bounded distance from a straight line.
	\item[(T2)] $\Omega$ cannot have only one planar strip end.
	\item[(T3)] If $\Omega$ has exactly two planar strip ends, then there exists a line $L$ such that $\Omega$ is at bounded distance from $L$, and the two ends are on opposite sides with respect to any line orthogonal to $L$.
\end{enumerate}
\end{theorem}

The similarity of the statements of our results on $f$-extremal domains and the parallel result in the context of constant mean curvature surfaces is evident. Nevertheless, the two problems are very different, we only recall the fact that the geometry of constant mean curvature surfaces is local, while the geometry of $f$-extremal domains is global! Consequently, the proofs of such two parallel results are different.

\medskip

%In the constant mean curvature case, W. H. Meeks in \cite{Meeks}, and N. J. Korevaar, R. Kusner and B. Solomon in \cite{KKS} proved that finite topoloy properly embedded surfaces have ends asymptotic to Delaunay surfaces. The construction of surfaces with this geometry has led to powerful methods in Geometric Analysis, see the paper of N. Kapouleas \cite{kap} and the papers of R. Mazzeo, F. Pacard and D. Pollack \cite{maz-pac, maz-pac-pol}. A particularly related situation to our overdetermined problem (\ref{pr_bis}) is the study of coplanar end surfaces, see the paper of K. Gro$\beta$e-Brauckmann, R. Kusner and J. Sullivan \cite{gro-kus-sul} and that of these authors joint with N. Korevaar and J. Ratzkin \cite{gro-kor-kus-rat-sul}.\\
%For $n=2$, Theorem \ref{T0} says that if the domain $\Omega$ is unbounded and contained in a strip, then it is symmetric with respect to an axis contained in the strip and its boundary has two connected components which are graphs over that axis. We conjecture that if $\Omega$ has proper finite topology, then each planar strip end of $\Omega$ has this asymptotic geometry.\\

Now, we want to link property $P_1$ with the function $f$ that appears in the overdetermined problem \eqref{pr_bis}. Consider a Lipschitz function $f$ that satisfies the property
\begin{enumerate}
\item [$P_2$] : there exists a positive constant $\lambda$ such that $f(t) \geq \lambda\,t$ for all $t>0$.
\end{enumerate}

Using basically the maximum principle and the moving plane argument we will obtain the following:

\begin{theorem}\label{T1}
Let $\Omega$ be an $f$-extremal domain of $\mathbb{R}^2$, where $f$ satisfies property $P_2$. Then, the following properties hold:
\begin{enumerate}
\item[(T4)] There exists a positive constant $R$ such that $\overline{\Omega}$ does not contain any closed ball of radius $R$ (i.e., $\Omega$ satisfies property $P_1$).
\item[(T5)] There exists a positive constant $h_0$ such that every connected component of
\[
\{ x \in \Omega\, |\, u(x) > h_0\},
\]
where $u$ is the solution of (\ref{pr_bis}), 
is contained in a ball of radius $\frac{\sqrt{5}}{2}\,R$.
\end{enumerate}
\end{theorem}

In fact, as we will prove, Statement (T4) is true for $f$-extremal domains $\Omega$ contained in $\mathbb{R}^n$, and not only in $\mathbb{R}^2$.
By Theorem \ref{T1}, property $P_1$ is satisfied if $f$ satisfies property $P_2$. Then all the conclusions of Theorem \ref{T2} are true for $f$-extremal domains of $\mathbb{R}^2$ when $f$ satisfies property $P_{2}$. In particular, we obtain the following result, that gives a partial affirmative answer to the conjecture of Berestycki-Caffarelli-Nirenberg in dimension 2:

\begin{corollary}\label{cor_BCN}
Let $\Omega$ be an $f$-extremal domain of $\mathbb{R}^{2}$, where $f$ satisfies property $P_{2}$, such that $\mathbb{R}^2 \backslash \overline{\Omega}$ is connected. Then, $\Omega$ is a ball.
\end{corollary}

Corollary \ref{cor_BCN} follows immediately from Theorem \ref{T2} because if
$\Omega$ is a (connected) domain and its complement is connected, then
$\Omega$ has the topology of a disc and its boundary consists just of one planar
curve. Hence $\pp \Omega$ separates $\mathbb{R}^{2}$ into two connected components. Then, $\Omega$ is either bounded, or the complement of a compact domain, or a proper finite topology domain with only one end. The corollary now follows from Theorems \ref{T1} (T4), \ref{T2} (T2), and the classical Serrin's result. We emphasize the fact that boundedness of the solution $u$ of the overdetermined problem \eqref{pr_bis} is not assumed: under the hypothesis of Corollary \ref{cor_BCN}, $u$ is bounded ``a fortiori''. Moreover, we remark that under the hypothesis $P_2$, the complement of a ball and the half-plane do not support any solution (bounded or unbounded) to the overdetermined problem \eqref{pr_bis}. %Another way of giving a partial answer to the conjecture of Berestycki-Caffarelli-Nirenberg, that is also a consequence of Theorems \ref{T1} (T4) and \ref{T2} (T2) is the following:
 
% \begin{corollary}
%Let $\Omega$ be an $f$-extremal unbounded domain of $\mathbb{R}^{2}$ such that $\mathbb{R}^2 \backslash \overline{\Omega}$ is connected. Then, $f$ does not satisfy property $P_{2}$ and $\Omega$ contains arbitrarily large balls.
%\end{corollary}

\medskip

Let come back to constant mean curvature (hyper-)surfaces. In $\mathbb{R}^3$, we saw that if the surface has exactly two ends, then it is contained in a cylinder (Theorem \ref{Meeks_T}). In this case, the geometry of the surface is very special. The second result about constant mean curvature (hyper-)surfaces we are interested in is the very well known:

\begin{theorem}\label{T_KKS} 
(N. J. Korevaar, R. Kusner and B. Solomon, \cite{KKS}). If $M$ is a properly embedded, non-zero constant mean curvature surface $M$ in $\mathbb{R}^3$ (or more generally a properly embedded, non-zero constant mean curvature hypersurface in $\mathbb{R}^n$) contained in a solid cylinder $B_R^n \times \mathbb{R}$ for some positive $R$, then $M$ is rotationally symmetric with respect to a line parallel to the axis of the cylinder $\{0\} \times \mathbb{R}$. 
\end{theorem}

In this paper we will show that a parallel result can be stated for $f$-extremal domains in $\mathbb{R}^n$, for $n \geq 2$, contained in solid cylinders (or in a planar strip if $n=2$). Very surprisingly, when $n = 2$ the hypothesis that the domain is contained in a planar strip can be replaced by the hypothesis that the domain is contained in a half-plane, and this will lead to a stronger result in dimension 2 that will allow us to give an other partial answer to the conjecture of Berestycki-Caffarelli-Nirenberg. In order to state the result, we give some definitions.

\begin{definition} Let $\Omega$ be a domain whose boundary is made by a unique proper arc $\Gamma$.\\ 
We say that the domain $\Omega$ is an \textbf{epigraph} if, after a suitable choice of coordinates, $\Gamma$ is the graph of a $C^2$ function $\varphi:\mathbb{R}\longrightarrow \mathbb{R}$, i.e., 
\[
\Omega =\{(x,y)\in \mathbb{R}^{2}\, |\, y>\varphi(x) \}.
\] 
The epigraph $\Omega$ is said to be \textbf{coercive} if 
\[
\lim_{|x|\rightarrow +\infty} \varphi(x) =+\infty.
\]
The epigraph $\Omega$ is said to be \textbf{uniformly Lipschitz} if $\varphi$ is uniformly Lipschitz.\\
We say that the domain $\Omega$ is an \textbf{almost epigraph} if $\Gamma = \{(x(t), y(t)) \, , \, t \in \mathbb{R}\}$ with $x'(t) \geq 0$ and $\{x(t) \, , \, t \in \mathbb{R}\} = \mathbb{R}$.
%We say that $\Omega$ is a \textbf{weak epigraph}, with respect to a straight line $L$, if the orthogonal projection of $\Gamma$ over $L$ is onto and monotonic.  
\end{definition}

%If $\gamma(t)=(x(t),y(t)), t\in\mathbb{R}$, is an arc-length parametrization of $\Gamma$,  $\gamma'=(x',y') the unit tangent vector and $n=-\nu=(-y',x')$ the inward pointing normal vector along  $\Gamma$, note that the assumption $x'>\varepsilon$ for some $\varepsilon >0$ implies that $\Omega$ is a uniformly Lipschitz epigraph. %Moreover, $\Omega$ is a weak epigraph with respect to $L = \{y=0\}$ if and only if the image of $x$ is the whole $\mathbb{R}$ and $x'\geq 0$. 

The result we will prove is the following:

\begin{theorem}\label{T4}Let $\Omega$ be an $f$-extremal domain of $\mathbb{R}^{2}$ (no restriction about the topology of the domain, and $f$ is only supposed to be a Lipschitz function). The following properties hold:
\begin{enumerate}
\item[(T6)] If $\Omega$ is contained in a wedge of angle less than $\pi$, then $\Omega$ is either a ball or a uniformly Lipschitz epigraph. 
\item[(T7)] If $\Omega$ is contained in a wedge of angle less than $\pi/2$, then $\Omega$ is a ball.
\item[(T8)] If $\Omega$ is contained in a half-plane, then $\Omega$ is either a ball or (after
a rigid motion) there exists a $C^2$ positive function
$\varphi:\mathbb{R}\longrightarrow ]0,\infty[$ such that either
\begin{enumerate}
\item[i.]the domain $\Omega$ is an epigraph $\{y>\varphi (x)\}$, or
\item[ii.] $\varphi$ is bounded and $\Omega$ is the symmetric domain
$\{|y|<\varphi(x)\}$.
\end{enumerate}
\item[(T9)] If $\Omega$ is unbounded and $\partial \Omega$ consists of a unique proper arc (we recall that an arc is proper if the intersection of the arc with any compact ball is compact), then either $\Omega$ is an almost epigraph or it contains a half-plane.
\end{enumerate}
\end{theorem}

In particular, note that under the hypothesis of Statement (T9) $f$ does not satisfy property $P_2$. Remark that Statement (T8) of the previous result is in some sense the parallel of the result of N. J. Korevaar, R. Kusner and B. Solomon for $f$-extremal domains in $\mathbb{R}^2$, and in fact it is much stronger because we suppose only that the domain lies in a half-plane and not in a planar strip. For the other dimensions we have the:

\begin{theorem} \label{T0} Let $\Omega$ be an $f$-extremal domain of $\mathbb{R}^n$, $n>2$, satisfying the property that there exists a line $L$ such that $\Omega$ is at bounded distance from $L$ (i.e., $\Omega$ is contained in a cylinder). Then, $\Omega$ has two ends, its boundary is rotationally symmetric with respect to a straight line parallel to $L$ and its generating curve is a bounded planar graph over this axis, i.e. there exists a $C^2$ positive function
$\varphi:\mathbb{R}\longrightarrow ]0,\infty[$ such that $\Omega$ (after a suitable rigid motion) is the domain
$\{(x,y) \in \mathbb{R}\times \mathbb{R}^{n-1}\,\,\, |\,\?\? |y|<\varphi(x)\}$.
\end{theorem}

The case of dimension 2 is in fact much more interesting and some remarks are due. Firstly, Corollary \ref{cor_BCN} can be obtained also from Statement (T9) of Theorem \ref{T4} using Statement (T4) of Theorem \ref{T1}, and Statement (T2) of Theorem \ref{T2} can be obtained also directly from Statement (T7) of Theorem \ref{T4}. Moreover, it is clear that if $\Omega$ is an epigraph contained in a wedge $\{y > c|x|\}, c>0$ (i.e., the angle of the wedge is less than $\pi$), then it is a coercive epigraph. In \cite{Far-Val-0}, A. Farina and E. Valdinoci proved that if $f$ is locally Lipschitz, $n=2$, and $\Omega$ is a uniformly Lipschitz coercive epigraph of class $C^3$, then there exists no function $u \in C^2(\overline{\Omega}) \cap L^\infty(\Omega)$ which is solution of (\ref{pr_bis}). Then, from Statement (T6) of Theorem \ref{T4} we have immediately the following:

 \begin{corollary}
 Let $\Omega$ be a 2-dimensional $C^3$-domain (with arbitrary topology) and $u$ a solution
 of problem (\ref{pr_bis}). If $\Omega$ is contained in
 a wedge of angle less than $\pi$ and $u$ is bounded, then $\Omega$ is
 a ball.
 \end{corollary}
 
%\textcolor{red}{An other interesting remark is that if $u$ is a solution of problem (\ref{pr_bis}), if $u$ is bounded then $|\nabla u|$ is bounded, see section 1.2 of \cite{Far-Val-0} and \cite{gilbarg}.} 
In our proof of Statement (T8) we use a tilted moving line argument and in the case of item (i) we obtain that the moving line reflection can be applied for any horizontal line. This fact shows that if $\Omega$ is a coercive epigraph, then $u$ must be strictly increasing in the second variable, see \cite{esteban} and \cite{BCN}. In \cite{Far-Val-0}  A. Farina and E. Valdinoci proved that if $f$ is locally Lipschitz, $n=2$, and $u \in C^2(\overline{\Omega})$ is a solution of (\ref{pr_bis}) increasing in the second variable and with bounded gradient, then $\Omega$ is a half-plane. From Statement (T8) of Theorem \ref{T4} we have immediately the following:

%\begin{remark} Even for unbounded $u$, it is natural to conjecture that any $\Omega$ contained in a wedge of angle less than $\pi$ is a ball. \end{remark}

%A basic geometric step of Theorem 1.6 of \cite{Far-Val-0} is the fact that if $\Omega$ is a coercive epigraph, then $u$ must be strictly increasing in the second variable. This follows from the moving line argument, see \cite{esteban} and \cite{BCN}.
%In our proof of T8 we use a tilted moving line argument and in the case of item i we obtain that the moving line reflection can be applied for any horizontal line. It follows that $\partial u/\partial y > 0$ and from Theorem 1.2 in \cite{Far-Val-0} we have the

\begin{corollary}\label{BCNhalfplane}
Let $\Omega$ be a $C^3$-domain in $\mathbb{R}^2$ such that $\mathbb{R}^2\backslash \Omega$ is
connected and $u$ be a solution of problem (\ref{pr_bis}). If $\Omega$ is contained in
a half-plane and the gradient of $u$ is bounded, then either $\Omega$ is
a ball and $u$ is a radial function or (after a rigid motion)
$\Omega =\{y>0\}$ and $u$ depends only on the variable $y$.
\end{corollary}

Corollary \ref{BCNhalfplane} proves the Berestycki-Caffarelli-Nirenberg conjecture when $\Omega$ lies in a half-plane and $|\nabla u|$ is bounded. Corollary \ref{cor_BCN} and Corollary \ref{BCNhalfplane} prove Theorem \ref{Tmain}.

\begin{remark}\label{nablau_acotado}
For an epigraph $\Omega$ there are two geometric situations where the condition that $u$ is bounded implies that $|\nabla u|$ is bounded: the case when the curvature of $\partial \Omega$ is bounded, and the case when $\Omega$ is a uniformly Lipschitz epigraph. In the first case, assuming that $\nabla u$ is unbounded, it is possible to choose a sequence of points $p_n \in \Omega$ with $|\nabla u(p_n)| \to +\infty$, such that the homotheties moving $p_n$ to the origin and rescaling $|\nabla u(p_n)|$ to 1 transform the function $u : \Omega \to \mathbb{R}$ into a sequence of functions $v_n : D_n \to \mathbb{R}$ with $|\nabla v_n(0)|=1$ and uniformly bounded gradient on compact subsets. Passing to a subsequence we get at the limit a bounded nonnegative function $v$ defined on a region $D \subset \mathbb{R}^2$ either without boundary or bounded by a straight line satisfying 
\[
\left\{\begin{array} {ll}
\Delta v = 0 & \mbox{in }\; D\\
      v= 0 \, \, \, \textnormal{and}\,\, \,  \langle \nabla u, \nu \rangle = 0& \mbox{on }\; \pp D \\
  |\nabla v| = 1 &  \mbox{at }\; 0 \in D\\
\end{array}\right.
\]
Liouville's theorem implies that $D \neq \mathbb{R}^2$. Then $D$ is a halfplane and in this case we contradict the unique continuation principle along $\partial D$.\\
In the case when $\Omega$ is a uniformly Lipschitz epigraph and $u$ is bounded, the boundedness of $\nabla u$ follows from the results of section 6.2 of \cite{gilbarg}.\footnote{Authors wish to thank Alberto Farina for pointing out to them this last fact.}
\end{remark}

\medskip

In order to present the third part of this paper, we come back once again to constant mean curvature (hyper-)surfaces. In \cite{Ros-Ros}, A. Ros and H. Rosenberg study some global properties of constant mean curvature surfaces in $\mathbb{R}^3$ that are contained in a slab (i.e. between two parallel planes). The natural class of surfaces to be considered in this case is the class $S$ of properly embedded non-zero constant mean curvature surfaces $M$ satisfying the following conditions: $M$ lies in the slab between two horizontal planes $\pi_1$ and $\pi_2$, is symmetric about the plane $\pi_0 = \{(x,y,z) \in \mathbb{R}^3\, |\, z=0\}$ and $M \cap \{z>0\}$ is a graph over the open subset $\Omega \subset \pi_0$. The infimum of the distance between $\pi_1$ and $\pi_2$ is called {\it width} of the surface $M$. The reason for which it is natural to consider such class of surfaces is well explained in  \cite{Ros-Ros}. The third result about constant mean curvature theory we are interested in is the following:

  \begin{theorem} \label{rosros} (A. Ros and H. Rosenberg, \cite{Ros-Ros}). Suppose $M \in S$ and let $H >0$ its mean curvature. If $M$ has width less than $1/H$, then the components of $\pi_0 \backslash \Omega$ are strictly convex. In particular $M$ is connected. If moreover $M$ has bounded curvature, then $\pi_0 \backslash \Omega$ is a countable disjoint union of strictly convex compact disks.
  \end{theorem}

In the last section of this paper, we will prove a kind of parallel of the previous result for double periodic domains in $\mathbb{R}^2$ where system \eqref{pr_bis} can be solved. Let
$\mathbb{T}^2$ be a flat torus obtained as a quotient of $\mathbb{R}^2$ by a
lattice, i.e., $\mathbb{T}^2 = \mathbb{R}^2 / \langle v_1, v_2 \rangle$,
where $v_1$ and $v_2$ are two linearly independent non-zero vectors of
$\mathbb{R}^2$ and

\[
\langle v_1,v_2\rangle = \{a\,v_1 + b\,v_2\,\, :\,\, a,b \in \mathbb{Z} \}.
\]

It is clear that a (connected) domain in $\mathbb{T}^2$ corresponds to a
(possibly nonconnected) double periodic domain in $\mathbb{R}^2$.
We will prove the following result:

\begin{theorem}\label{T3}
Let $\Omega$ be a (connected) domain of $\mathbb{T}^2$ that supports a
solution $u \in C^2(\overline{\Omega}) \cap C^3(\Omega)$ to problem
(\ref{pr_bis}), where $f$ is a $C^1$ function such that
\begin{equation}\label{hyp}
2\, \max_{x \in\overline{\Omega}}\,\int_0^{u(x)}f(s)\,\textnormal{d}s %\leq
< \alpha^2.
\end{equation}
Then, each component of $\mathbb{T}^2 \backslash \Omega$ is strictly convex.
%or $\Omega$ is the quotient in  the torus $\mathbb{T}^2$ of a strip of the plane $\mathbb{R}^2$.
\end{theorem}

The previous theorem is true in particular for the linear function $f(t)=\lambda\,t$. We state it in the formulation of a double periodic domain of $\mathbb{R}^2$.

\begin{corollary}\label{cor_periodic}
Let $\Omega$ be a connected component of a (possibly nonconnected) double periodic open domain of $\mathbb{R}^2$. Let us suppose that there exists a doubly periodic solution $u \in C^2(\overline{\Omega}) \cap C^3(\Omega)$ to the overdetermined problem
\begin{equation}\label{pr_2}
\left\{
\begin{array} {ll}
\Delta\, u + \gl\, u = 0 &\mbox{in }\; \Omega\\
u > 0  & \mbox{in }\; \Omega \\
u=0 & \mbox{on }\; \pp \Omega \\
\langle \nabla u, \nu \rangle = \alpha &\mbox{on }\; \pp \Omega
\end{array} 
\right.
\end{equation}
where $\lambda$ is a positive constant and $\alpha$ is a constant. If
\begin{equation}
\max_{\Omega} u %\leq 
< \frac{|\alpha|}{\sqrt{\lambda}},
\end{equation}
then each component of $\mathbb{R}^2 \backslash \Omega$ is strictly convex. %or $\Omega$ is a strip.
\end{corollary}

We state explicitly Corollary \ref{cor_periodic} because it is interesting to remark that $\frac{|\alpha|}{\sqrt{\lambda}}$ is the maximum value of the function $u$ that satisfies \eqref{pr_2} in the strip $[0,\pi/\sqrt{\lambda}] \times \mathbb{R}$. Such a result implies that the maximum value of the (bounded) function $u$ that satisfies \eqref{pr_2} in the domains of the family $s \to \Omega_{s}$ of periodic and symmetric perturbations of the strip found by F. Schlenk and P. Sicbaldi in \cite{Sch-Sic}, is bigger or equal to $\frac{|\alpha|}{\sqrt{\lambda}}$. %and equal to $\frac{|\alpha|}{\sqrt{\lambda}}$ only when $\Omega_{s} = \Omega_{0}$, i.e. the strip. 
This follows from the fact that $\mathbb{R}^2 \backslash \Omega_s$ is not convex (by the construction of $\Omega_s$).

\b
\ni
{\bf Digressions, new ideas and open problems.} With the previous results it is clear that the link between ovedetermined elliptic problems and constant mean curvature surfaces if very strong.
In the constant mean curvature case, W. H. Meeks in \cite{Meeks}, and N. J. Korevaar, R. Kusner and B. Solomon in \cite{KKS} proved that the only unbounded constant mean curvature surfaces contained in a cylinder are the Delaunay surfaces, and that finite topology properly embedded surfaces have ends asymptotic to Delaunay surfaces. The development of the technique of surface gluing and the construction of constant mean curvature surfaces with Delaunay ends has led to powerful methods in Geometric Analysis, see the paper of N. Kapouleas \cite{kap} and the papers of R. Mazzeo, F. Pacard and D. Pollack \cite{maz-pac, maz-pac-pol}. A related situation is the study of coplanar end surfaces, see the papers by C. Cosin and A. Ros \cite{cos-ros}, K. Gro$\ss$e-Brauckmann, R. Kusner and J. Sullivan \cite{gro-kus-sul} and that of these authors joint with N. Korevaar and J. Ratzkin \cite{gro-kor-kus-rat-sul}. An other very powerful tool in Geometric Analysis is the use (often for comparison with the maximum principle) of the catenoid, a very well known unbounded minimal surface of revolution. It is important to remark that the family of Delaunay surfaces can be parameterized by $\sigma \in (0,1]$, where the value $\sigma= 1$ corresponds to the straight cylinder, and the limit value $\sigma \to 0$ corresponds to the union of spheres. In this situation the catenoid can be seen as the limit for $\sigma \to 0$, in some sense to be made precise, of the neck of the Delaunay surfaces. With respect to overdetermined elliptic problems, it is tempting to conjecture that:
\begin{itemize}
\item The family of Delaunay type extremal domain $s \to \Omega_s$ constructed by F. Schlenk and P. Sicbaldi in \cite{Sch-Sic} for $s$ in a small interval $(1-\epsilon, 1]$ can be in fact parametrized by the interval $s \in (0,1]$ and the limit case $s \to 0$ corresponds to the union of balls.
\item There exists an unbounded $0$-extremal domain that looks like the region inside of a catenoid. This idea comes from the fact that $0$-extremal domains arise as limits under scaling of sequences of $\lambda$-extremal domains, and for this last problem we have the starting point of the Delaunay type family of extremal domains $\{\Omega_s\}_{s \in (1-\epsilon,1]}$. In $\mathbb{R}^2$ such domain is $\Omega_*$ defined in the previous section, but the question is open in $\mathbb{R}^n$ for $n\geq 3$. Remark that the catenoid in $\mathbb{R}^n$, $n \geq 4$ is contained in a slab. Then, by the shift of dimensions that occurs when one consider extremal domains instead of constant mean curvature hypersurfaces, the boundary of the $0$-extremal domain of $\mathbb{R}^n$ that looks like the region inside of a catenoid should be contained in a slab for $n\geq 3$. It is interesting to remark also that in \cite{nuovo} D. Khavinson, E. Lundberg and R. Teodorescu prove that there does not exist a rotationally-symmetric 0-extremal domain in $\mathbb{R}^4$ that contains its own axis of symmetry and whose boundary is obtained by rotating the (two-dimensional) graph of an even real-analytic function about the $x$-axis.
\item In the case of finite topology, the ends of an $f$-extremal domain with $f(t) = \lambda\, t$ are asymptotic to the ends of the expected family of Delaunay type domains $\{\Omega_s\}_{s \in (0,1]}$, and a similar geometry can be obtained also for $f$-extremal domains with $f(t) \geq \lambda\, t$.
\item It is possible to construct highly non-trivial $f$-extremal domain of finite topology with ends asymptotic to Delaunay type domains $\Omega_s$.
\end{itemize} 
As a final remark, we observe that the real question to solve is the characterization of overdetermined elliptic problems, i.e. to find necessary and sufficient conditions on the geometry of a general $f$-extremal domain, and to find necessary and sufficient conditions on the function $f$ and the constant $\alpha$ to get existence of a solution to the elliptic system (\ref{pr_bis}). This question is highly non-trivial. This paper can help in the following sense: with our results, one has necessary conditions on the geometry and the topology of general $f$-extremal domains, and in fact there are few possibilities of such domains. Then, it should be not so hard to check, for each one of such possible domains, which are the good functions $f$ and the good constants $\alpha$ that give existence of a solution of (\ref{pr_bis}).

%A particularly related situation to our overdetermined problem (\ref{pr_bis}) is the study of coplanar end surfaces, see the papers of C. Cosin and A. Ros \cite{cos-ros}, and K. Gro$\beta$e-Brauckmann, R. Kusner and J. Sullivan \cite{gro-kus-sul} and that of these authors joint with N. Korevaar and J. Ratzkin \cite{gro-kor-kus-rat-sul}. From property T8 of theorem \ref{T4} we have that if the domain $\Omega$ is unbounded and contained in a strip, then it is symmetric with respect to an axis contained in the strip and its boundary has two connected components which are graphs over that axis. We conjecture that if $\Omega$ has proper finite topology, then each planar strip end of $\Omega$ has this asymptotic geometry. We remark that the domains of the family $s \to \Omega_{s}$ of periodic and symmetric perturbations of the strip found by F. Schlenk and P. Sicbaldi in \cite{Sch-Sic} are examples of such kind of extremal domains.

\b
\ni
{\bf Plan of the paper.} In order to simplify the exposition, we will start by proving Theorem \ref{T1}, then we will continue with Theorems \ref{T2}, \ref{T4}, \ref{T0}, and \ref{T3}.

\section{A narrowness property of the domain via the maximum principle}

To every positive constant $\lambda$ we can associate the radius $R_{\lambda}$ of balls whose first eigenvalue of the Laplacian with 0 Dirichlet boundary condition is $\lambda$. In other words, $R_{\lambda}$ is the positive constant such that one can solve
\begin{equation}\label{pr:ball}
\left\{
\begin{array} {ll}
\Delta\, v + \gl\, v = 0 &\mbox{in }\; B_{R_{\lambda}}(p)\\
v > 0  & \mbox{in }\; B_{R_{\lambda}}(p) \\
v=0 & \mbox{on }\; \pp B_{R_{\lambda}}(p)
\end{array} 
\right.
\end{equation}
where $B_{R_{\lambda}}(p)$ is the ball of $\mathbb{R}^n$ of radius ${R_{\lambda}}$ and center $p \in \mathbb{R}^2$. 
 We remark that the first eigenvalue of the Dirichlet-Laplacian on a ball
of radius $R$ is given by
\[
R^{-2}\, \lambda_1(B_1)
\]
where $\lambda_1(B_1)$ is the first eigenvalue the Laplacian on the unit
ball in $\mathbb{R}^n$ with 0 Dirichlet boundary condition, and then
$R_{\lambda}$ depends on the constant $\lambda$ and the dimension $n$.

\medskip

Statement (T4) of Theorem \ref{T1}, and some more details, are an immediate
consequence of the following result:

\begin{proposition}\label{balls_R}
Let us suppose $\Omega$ is an open (bounded or unbounded) connected domain of $\mathbb{R}^n$ such that one can find a (strictly) positive function $u \in C^2(\Omega)$ that solves the elliptic equation
\begin{equation}\label{elliptic}
\begin{array} {ll}
\Delta u + f(u) = 0 & \mbox{in }\; \Omega,\\
						\end{array}
\end{equation}
where $f : (0,+\infty) \to \mathbb{R}$ satisfies property $P_2$. Then, $\Omega$ does not contain any closed ball of radius $R_{\lambda}$. Moreover, if $u$ satisfies the boundary conditions
\begin{equation}\label{boundary}
\begin{array} {ll}
 u = 0 & \mbox{on }\; \pp \Omega\\
 \langle \nabla u, \nu \rangle = \alpha &\mbox{on }\; \pp \Omega 
						\end{array}
\end{equation}
for some negative constant $\alpha$, then either the closure $\overline \Omega$ does not contain any closed ball of radius $R_{\lambda}$ or $\Omega$ is a ball of radius $R_{\lambda}$.
\end{proposition}

\textit{Proof.} Let $u$ be a solution of equation (\ref{elliptic}) where $f$ satisfies property $P_2$. Let us suppose that there exists a point $p \in \mathbb{R}^2$ such that $\overline{B_{R_{\lambda}}(p)} \subseteq \Omega$.  If $v$ is the solution of (\ref{pr:ball}) normalized to have $L^{2}$-norm equal to 1, then it is possible to choose $\epsilon > 0$ such that the function
\[
v_\epsilon = \epsilon\, v
\]
has the following properties (see fig. 1):
\begin{enumerate}
	\item $v_{\epsilon}(x) \leq u(x)$ for all $x \in \overline{B_{R_{\lambda}}(p)}$;
	\item there exists $x_{0} \in B_{R_{\lambda}}(p)$ such that $v_{\epsilon}(x_{0}) = u(x_{0})$.
\end{enumerate}
The function $v_\epsilon$ satisfies (\ref{pr:ball}) and then, by property $P_2$, the function $u-v_{\epsilon}$ satisfies
\begin{equation} 
\begin{array} {ll}
\Delta (u-v_{\epsilon}) \leq - \lambda\,(u-v_{\epsilon}) \leq 0 & \mbox{in }\; \overline{B_{R_{\lambda}}(p)}\\
						\end{array}
\end{equation}
Moreover  $u-v_{\epsilon}$ is a nonnegative function that attains its minimum at the interior point $x_{0}$. The maximum principle (see \cite{gilbarg}, p.32) leads to a contradiction.\\
\makeatletter
\begin{figure}[!ht]
\centering
{\scalebox{.8}{%\input{image1.pstex_t}
\begin{picture}(0,0)%
\includegraphics{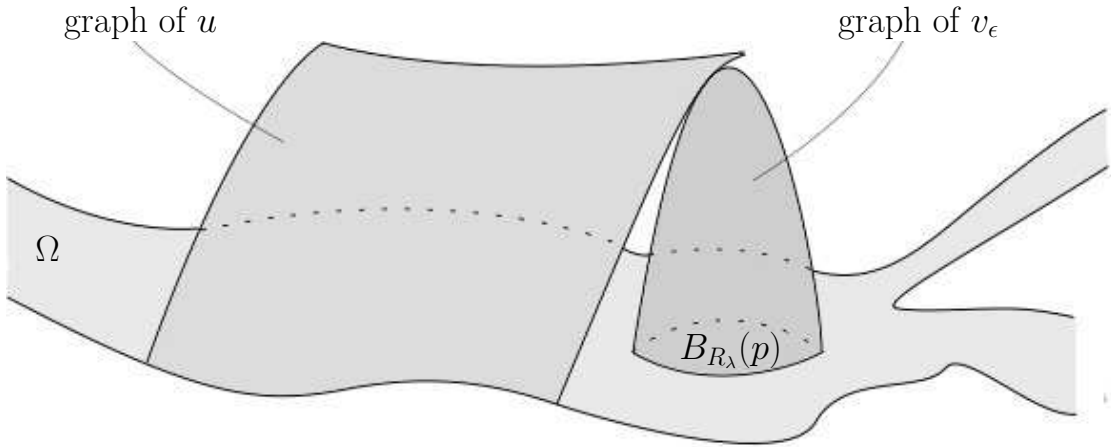}%
\end{picture}%
\setlength{\unitlength}{3947sp}%
\begingroup\makeatletter\ifx\SetFigFont\undefined%
\gdef\SetFigFont#1#2#3#4#5{%
  \reset@font\fontsize{#1}{#2pt}%
  \fontfamily{#3}\fontseries{#4}\fontshape{#5}%
  \selectfont}%
\fi\endgroup%
\begin{picture}(8730,3546)(751,-3001)
\put(976,-1486){\makebox(0,0)[lb]{\smash{{\SetFigFont{17}{20.4}{\familydefault}{\mddefault}{\updefault}{\color[rgb]{0,0,0}$\Omega$}%
}}}}
\put(6031,-2260){\makebox(0,0)[lb]{\smash{{\SetFigFont{17}{20.4}{\familydefault}{\mddefault}{\updefault}{\color[rgb]{0,0,0}$B_{R_\lambda}(p)$}%
}}}}
\put(1201,314){\makebox(0,0)[lb]{\smash{{\SetFigFont{17}{20.4}{\familydefault}{\mddefault}{\updefault}{\color[rgb]{0,0,0}graph of $u$}%
}}}}
\put(7276,314){\makebox(0,0)[lb]{\smash{{\SetFigFont{17}{20.4}{\familydefault}{\mddefault}{\updefault}{\color[rgb]{0,0,0}graph of $v_\epsilon$}%
}}}}
\end{picture}%
}}
\label{fig1}\caption{The domain $\Omega$, a portion of the graph of the function $u$ and the graph of the function $\epsilon\, v$ over the ball of radius $R_{\lambda}$.}
\end{figure}
\makeatother
Now let us suppose that $\Omega$ is unbounded and $u$ satisfies the boundary conditions (\ref{boundary}). By the previous statement it is clear that $\Omega$ does not contain any closed ball of radius $R_{\lambda}$. Let us suppose that there exists a point $p \in \mathbb{R}^2$ such that $\overline{B_{R_{\lambda}}(p)} \subseteq \overline \Omega$. Then, the boundary of $\overline{B_{R_{\lambda}}(p)}$ touches the boundary of $\Omega$ in some point $q$. Boundary conditions (\ref{boundary}) imply that there exists a positive constant $\delta_{0}$ such that the function 
\[
v_{\delta_{0}} = \delta_{0}\, v
\]
has the following properties:
\begin{enumerate}
	\item $v_{\delta_{0}}(x)< u(x)$ for all $x \in B_{R_{\lambda}}(p)$, and
	\item the Neumann data of $v_{\delta_{0}}$ at the boundary $B_{R_{\lambda}}(p)$ are equal to a constant $\beta$ such that $\alpha< \beta<0$.
\end{enumerate}
Now, as the parameter $\delta$ increases starting from $\delta_{0}$, defining $v_{\delta}$ as $\delta\,v$, one of the two situations occurs:
\begin{enumerate}
	\item $v_{\delta}(x_{0}) = u(x_{0})$ for some $x_{0} \in B_{R_{\lambda}}(p)$;	
	\item the Neumann data of $v_{\delta}$ becomes equal to $\alpha$ and $v_{\delta}(x)< u(x)$ for all $x \in B_{R_{\lambda}}(p)$.
\end{enumerate}
The first situation above implies, by the maximum principle, that $u = v_\delta$ and then $\Omega = B_{R_{\lambda}}(p)$. In the second case, we have that $u-v_{\delta}$ is a positive function in 
$B_{R_{\lambda}}(p)$ with
\[
\Delta (u-v_{\delta}) \leq - \lambda\,(u-v_{\delta}) \leq 0
\]
and at $q \in \pp \Omega \cup \pp B_{R_{\lambda}}(p)$ we have $(u-v_{\delta})(q) = 0$ and $\langle \nabla (u-v_{\delta}), \nu \rangle (q) = 0$, leading to a contradiction by the maximum principle (see \cite{gilbarg}, p.34). \proofend

The previous proposition says us that if $f$ satisfies property $P_2$, then the domain $\Omega$ is quite narrow, in the sense that it does not contain any ball which radius is bigger or equal to the given constant $R_{\lambda}$. An immediate consequence is the following:

\begin{remark}\label{imposs}
If $\Omega$ admits a positive solution of (\ref{elliptic}) and
$f$ satisfies property $P_2$, then $\Omega$ cannot be neither the complement of a
ball, nor a half-space, nor an epigraph, nor the complement of a cylinder
$B^k \times \mathbb{R}^{n-k}$ where $B^k$ is a round ball in
$\mathbb{R}^k$.

\end{remark}

\section{Symmetry properties of the domain via the moving plane method}

One of the most important tools coming from the maximum principle is the moving plane method. It was introduced by A. D. Alexandrov \cite{Alex} in order to prove that the only embedded, compact mean curvature hypersurface in $\mathbb{R}^n$ is the sphere. In a very elegant paper \cite{Serrin}, J. Serrin adapted the moving plane method to bounded domains where the elliptic overdetermined problem (\ref{pr_bis}) can be solved, in order to prove a strong symmetry property. In fact, he improved one of the central ingredients of Alexandrov's proof, the maximum principle at the boundary, proving what we now call the boundary maximum principle at a corner. Let us outline the result of J. Serrin \cite{Serrin}, see also \cite{Puc-Ser}.\\
Let us suppose that $\Omega$ is a bounded open domain of $\mathbb{R}^n$ whose boundary is of class $C^2$ and there exists a solution $u \in C^2(\overline{\Omega})$ to problem (\ref{pr_bis}), where $f$ is of class $C^1$ (in fact only Lipschitz regularity is required, as shown in \cite{Puc-Ser}). Let $T_0$ be a hyperplane in $\mathbb{R}^n$ not intersecting the domain $\Omega$ (the boundedness of $\Omega$ guarantees the existence of $T_0$). We suppose this hyperplane to be continuously moved normal to itself until it intersects by first time $\Omega$. From that moment onward, at each stage of the motion the resulting hyperplane $T$ will cut off from $\Omega$ a bounded cap $\Sigma(T)$ ($\Sigma(T)$ is the portion of $\Omega$ which lies on the same side of $T$ as the original hyperplane $T_0$, and its boundedness comes from the boundedness of $\Omega$). For any cap $\Sigma(T)$ thus formed, let $\Sigma'(T)$ be its reflection about $T$. $\Sigma'(T)$ is contained in $\Omega$ at the beginning of the process, and indeed as $T$ advances into $\Omega$, the resulting cap $\Sigma'(T)$ will stay within $\Omega$ at least until one of the following two events occurs:
\begin{enumerate}
\item $\Sigma'(T)$ becomes internally tangent to the boundary of $\Omega$ at some point not on $T$, or
\item $T$ reaches a position where it is orthogonal to the boundary of $\Omega$ at some point. 
\end{enumerate}  
Denote the hyperplane $T$ when it reaches either one of these positions by $T'$. The main result of J. Serrin is the following: 

\begin{theorem}\textnormal{(\textbf{J. Serrin, 1971}, \cite{Serrin})} \label{serrin_result} The reflected cap $\Sigma'(T')$ coincides with the part of $\Omega$ on the same side of $T'$ as $\Sigma'(T')$; that is, $\Omega$ is symmetric about $T'$.
\end{theorem}

As a corollary of this theorem we have that the only bounded domains $\Omega$ where one can solve (\ref{pr_bis}) are balls. In fact, the boundedness of $\Omega$ implies that for any given direction of $\mathbb{R}^n$, there exists an hyperplane $T'$ normal to that direction such that $\Omega$ is symmetric about $T'$. Moreover, the construction of $\Omega$ as union of caps $\Sigma(T')$ and $\Sigma'(T')$ implies that $\Omega$ is simply connected. The only simply connected domains which have this symmetry property are the balls.\\
We will refer to Theorem \ref{serrin_result} as the Serrin's reflection method. We remark that the boundedness of $\Omega$ is used only to guarantee the existence of the original non-intersecting plane $T_0$ and the boundedness of the cap $\Sigma(T)$. We remark also that the regularity hypothesis on the boundary of $\Omega$ (it is asked to be of class $C^2$) is a technical hypothesis used in the proof of Theorem \ref{serrin_result}.

\medskip

We can use the Serrin's technique to obtain some symmetry results for unbounded domains. Let $\Omega$ be an \underline{unbounded} open domain of $\mathbb{R}^n$ whose boundary is of class $C^2$ and let $u$ be a $C^2(\overline{\Omega})$-solution to problem (\ref{pr_bis}). Let $L$ be a hyperplane in $\mathbb{R}^n$ that intersects $\Omega$, and let $L^+$ and $L^-$ be the two connected components of $\mathbb{R}^n \backslash L$. We are interested in the geometry of \underline{bounded} connected components of $\Omega \cap L^+$ or $\Omega \cap L^-$.
 
\begin{proposition}\label{compact}
Let us suppose that $\Omega \cap L^+$ has a bounded connected component $C$. Then, the closure of $\pp C \cap L^+$ is a graph over $\pp C \cap L$ (see fig. 2).
\end{proposition}

\begin{figure}[!ht]
\centering
{\scalebox{.8}{%\input{image2.pstex_t}
\begin{picture}(0,0)%
\includegraphics{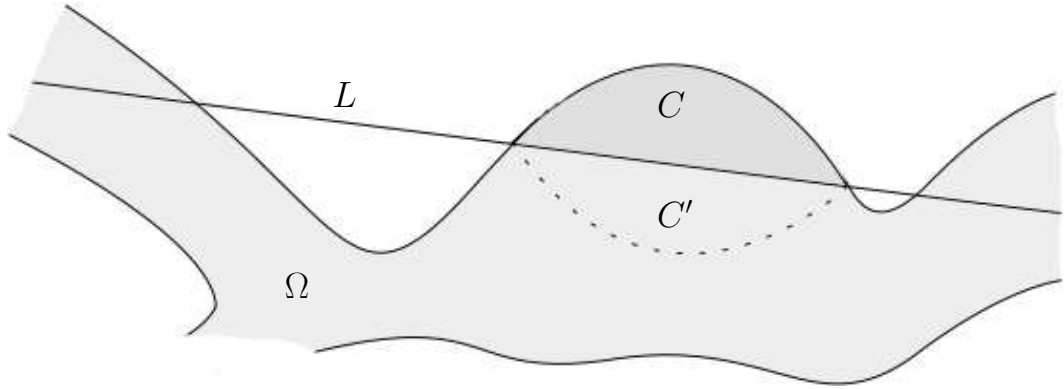}%
\end{picture}%
\setlength{\unitlength}{3947sp}%
\begingroup\makeatletter\ifx\SetFigFont\undefined%
\gdef\SetFigFont#1#2#3#4#5{%
  \reset@font\fontsize{#1}{#2pt}%
  \fontfamily{#3}\fontseries{#4}\fontshape{#5}%
  \selectfont}%
\fi\endgroup%
\begin{picture}(8340,3015)(901,-2701)
\put(6001,-586){\makebox(0,0)[lb]{\smash{{\SetFigFont{17}{20.4}{\familydefault}{\mddefault}{\updefault}{\color[rgb]{0,0,0}$C$}%
}}}}
\put(3076,-2011){\makebox(0,0)[lb]{\smash{{\SetFigFont{17}{20.4}{\familydefault}{\mddefault}{\updefault}{\color[rgb]{0,0,0}$\Omega$}%
}}}}
\put(3451,-511){\makebox(0,0)[lb]{\smash{{\SetFigFont{17}{20.4}{\familydefault}{\mddefault}{\updefault}{\color[rgb]{0,0,0}$L$}%
}}}}
\put(6001,-1456){\makebox(0,0)[lb]{\smash{{\SetFigFont{17}{20.4}{\familydefault}{\mddefault}{\updefault}{\color[rgb]{0,0,0}$C'$}%
}}}}
\end{picture}%
}}
\label{fig2}\caption{Moving plane method applied to the bounded component $C$.}
\end{figure}

\textit{Proof.} The proof of this proposition is based on the Serrin's reflection method. %Notice $L_t$ the 1-parameter family of straight hyperplanes in $L^+$, parallel to $L$, which distance from $L$ is equal to $t$, and notice $L^+_t$ and $L^-_t$ the two connected components of $\mathbb{R}^n \backslash L_t$, such that $L^+_t \subseteq L^+ = L^+_0$. 
By the boundedness of $C$, there exists a hyperplane $T_0 \in L^+$ parallel to $L$ not intersecting $C$. When this hyperplane is continuously moved normal to itself, it will intersect $C$ a first time. From that moment on, at each stage of the motion the resulting hyperplane $T$ will cut off from $C$ a bounded cap $\Sigma(T)$. For any cap $\Sigma(T)$ thus formed, let $\Sigma'(T)$ be its reflection about $T$.
$\Sigma'(T)$ is contained in $\Omega$ at the beginning of the process, and as $T$ advances into $C$, the resulting cap $\Sigma'(T)$ will stay within $\Omega$ at least until one of the following three events occurs:
\begin{enumerate}
\item $\Sigma'(T)$ becomes internally tangent to the boundary of $\Omega$ at some point not on $T$, or
\item $T$ reaches a position where it is orthogonal to the boundary of $\Omega$ at some point, or 
\item $T$ coincides with $L$.
\end{enumerate}  
The first two events are not possible by the Serrin's reflection, because $\Omega$ is unbounded. This means that $\Sigma'(T)$  stays within $\Omega$ for all hyperplane parallel to $T_0$ staying between $T_0$ and $L$ (fig. 2), and then the closure of $\pp \Sigma(L) \cap L^+$, i.e., the closure of $\pp C \cap L^+$ is a graph over $\pp C \cap L$. 
%is tangent to $\pp C \cap L^+$ and such that $C \subseteq L^-_{t_0}$. Moreover there exists a $\epsilon_0 > 0$ such that the reflection of $C \cap L^+_{t_0-\epsilon}$ with respect to $L_{t_0-\epsilon}$ is contained in $C$ for all $\epsilon \in [0,\epsilon_0]$ (this follows from regularity of $\Omega$). Now, when $\epsilon$ increases from $\epsilon_0$ to $t_0$, there are two possibility:
%\begin{enumerate}
%	\item The reflection of $C \cap L^+_{t_0-\epsilon}$ with respect to $L_{t_0-\epsilon}$ is contained in $\Omega$ for all $\epsilon \in [0,t_0]$. In this case it is clear that $\pp C \cap L^+$ is a graph over $\pp C \cap L$.
%	\item There exists a $\epsilon \in [\epsilon_0,t_0]$ such that the reflection of $C \cap L^+_{t_0-\epsilon}$ with respect to $L_{t_0-\epsilon}$ touches $\pp \Omega \cap L^-_{t_0-\epsilon}$ in an interior or boundary point. By the Serrin method (\cite{Serrin}, p.305--311) this means that $\Omega = C_{t_0-\epsilon} \cup C^*_{t_0-\epsilon}$, where $C_{t_0-\epsilon} = C \cap L^+_{t_0-\epsilon}$ and $C^*_{t_0-\epsilon}$ is its reflection with respect to $L_{t_0-\epsilon}$, leading to a contradiction with the non boundedness of $\Omega$.
%\end{enumerate}
%Only the first case is possible and the lemma follows at once.
\proofend

The previous proposition and its proof immediately imply the following properties:
\begin{corollary}
$\pp C \cap L$ is connected. 
\end{corollary}

\begin{corollary}\label{cor_orth}
The closure of $\pp C \cap L^+$ is not orthogonal to $L$ at any point. 
\end{corollary}

\textit{Proof.} In fact, if the closure of $\partial C \cap L^{+}$ meets $L$ orthogonally, then $\Omega$ is symmetric with respect to $L$, which contradicts the fact that $\Omega$ is unbounded.

\begin{corollary}\label{cor_union}
If $C'$ is the reflection of $C$ about $L$, then the closure of $C \cup C'$ stays within $\overline{\Omega}$, see fig. 2. 
\end{corollary}

We remark that $L$ is an arbitrary hyperplane such that there exists a bounded connected component $C$ of $\Omega \cap L^+$.

Now, let us suppose that $f$ satisfies property $P_2$. Then by Proposition \ref{balls_R} and Corollary \ref{cor_union} we have the following:

\begin{corollary}\label{halfAlex}
If $f$ satisfies property $P_2$, then it is not possible to construct a half-ball of radius $R_{\lambda}$ having base on $\pp C \cap L$ and staying within $C$.
\end{corollary}

Corollary \ref{halfAlex} follows immediately from the fact that  if $f$ satisfies property $P_2$ and $C'$ is the reflection of $C$ about $L$, then the closure of $C \cup C'$  cannot contain any closed ball of radius $R_{\lambda}$.

\section{Boundedness properties for the solution of the elliptic problem}

Let $\Omega$ be an open unbounded connected domain of $\mathbb{R}^2$
whose boundary is of class $C^2$, and such that there exists a function
$u \in C^2(\overline{\Omega})$ that solves the elliptic problem
(\ref{pr_bis}).
For the moment we suppose $\alpha \neq 0$.
Let $R_{\lambda}$ be the radius of the ball whose first eigenvalue of the
Dirichlet-Laplacian is $\lambda$, and $v$ a solution of
(\ref{pr:ball}) such that
$$
\langle \nabla v, \nu \rangle = \alpha
$$
at the boundary of $B_{R_{\lambda}} (p)$. Denote
$$
h_{0} := h_0(\lambda, \alpha) := \max_{B_{R_{\lambda}} (p)} v = v(p)
$$

Statement (T5) of Theorem \ref{T1} follows from the following proposition and remark. Similar geometric ideas were used by J. M. Espinar, J. A. G\'alvez and H. Rosenberg in \cite{esp-gal-rosen} in the context of constant curvature surfaces.

\begin{proposition}\label{prop-h0}

Let $f$ satisfy property $P_2$ and $\alpha \neq 0$. Let $\Omega'$ be a connected component of
\[
\{ x \in \Omega\, |\, u(x) > h_0\}
\]
Then, the diameter of ${\Omega'}$ is smaller than $2R_\lambda$.
In particular, there exists a point $p$ such that
$\overline{\Omega'} \subset B_{R}(p)$, where
$R=\frac{\sqrt{5}}{2}R_\lambda$.
\end{proposition}

\textit{Proof.}
%Let $h \geq h_{0}$ be a regular value of $u$ (by the Sard's theorem, almost all $h \geq h_0$ are regular values). Consider  Let $C$ be a connected component of $\pp\Omega'$. $C$ could be bounded or unbounded, see fig. 3.\\
First let us suppose that $\Omega'$ is bounded. Let
$d$ be its diameter, and suppose $d \geq 2 R_{\lambda}$. Let $q_1$ and $q_2$ be two points of $\overline{\Omega'}$ such that the distance between
$q_{1}$ and $q_{2}$ is bigger or equal to $2 R_{\lambda}$, and $C$ a curve in $\overline{\Omega'}$ joining $q_1$ and $q_2$, see fig. 3 (if $\Omega'$ is regular, $C$ can be taken in its boundary). Let $m$ be the mid point of the segment $\overline{q_{1}q_{2}}$, denote by $L_{1}$ the line containing
$\overline{q_{1}q_{2}}$ and by $L_{2}$ the line orthogonal to the segment
$\overline{q_{1}q_{2}}$ passing through $m$. Let 
$\Gamma= (L_{1} \backslash \overline{q_{1}q_{2}}) \cup C$ and denote by $H_{1}$ and $H_2$ the two connected components of $\mathbb{R}^{2} \backslash \Gamma$. Let $\Omega_{1}  = \Omega \cap H_{1}$.
Let $p \in L_{2} \cap H_{2}$ a point very far from $\Omega_{1}$ and consider the graph $G$ of the function $v$ defined on $B_{R_{\lambda}} (p)$ by (\ref{pr:ball}). Now let us translate the point $p$ along
the line $L_{2}$ in order to approach the domain $\Omega_{1}$.

\begin{figure}[!ht]
\centering
{\scalebox{.8}{
\begin{picture}(0,0)%
\includegraphics{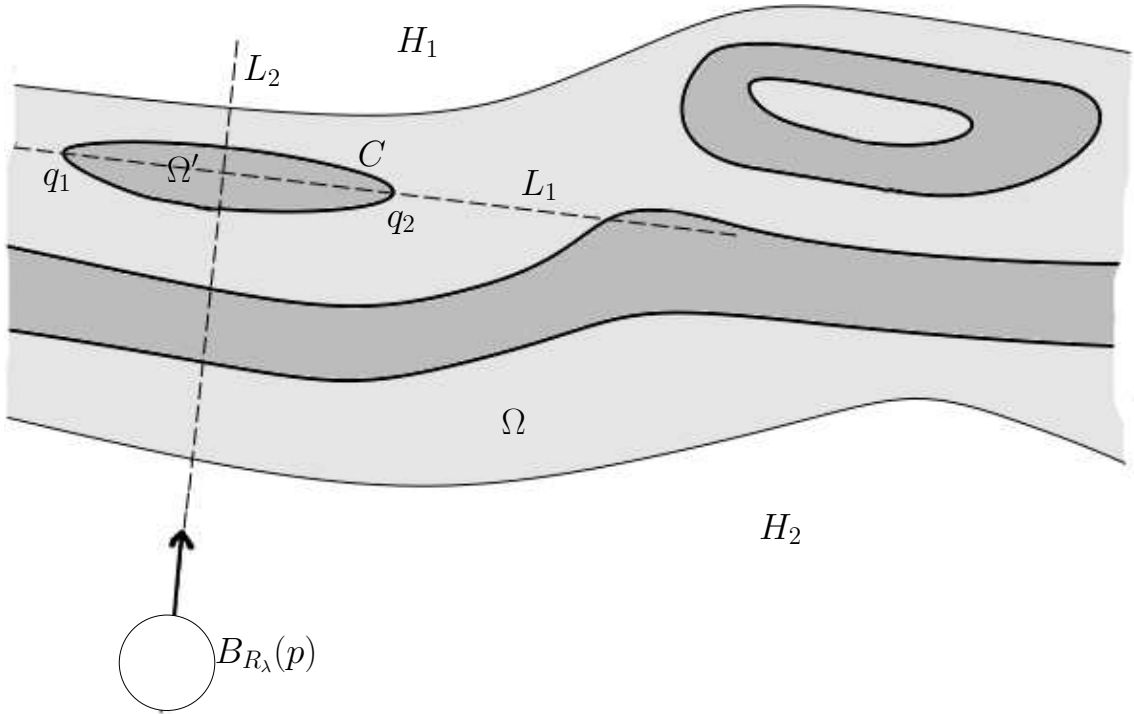}%
\end{picture}%
\setlength{\unitlength}{3947sp}%
\begingroup\makeatletter\ifx\SetFigFont\undefined%
\gdef\SetFigFont#1#2#3#4#5{%
  \reset@font\fontsize{#1}{#2pt}%
  \fontfamily{#3}\fontseries{#4}\fontshape{#5}%
  \selectfont}%
\fi\endgroup%
\begin{picture}(8910,5580)(526,-5416)
\put(826,-1261){\makebox(0,0)[lb]{\smash{{\SetFigFont{17}{20.4}{\familydefault}{\mddefault}{\updefault}{\color[rgb]{0,0,0}$q_1$}%
}}}}
\put(2401,-436){\makebox(0,0)[lb]{\smash{{\SetFigFont{17}{20.4}{\familydefault}{\mddefault}{\updefault}{\color[rgb]{0,0,0}$L_2$}%
}}}}
\put(4576,-1336){\makebox(0,0)[lb]{\smash{{\SetFigFont{17}{20.4}{\familydefault}{\mddefault}{\updefault}{\color[rgb]{0,0,0}$L_1$}%
}}}}
\put(3301,-1111){\makebox(0,0)[lb]{\smash{{\SetFigFont{17}{20.4}{\familydefault}{\mddefault}{\updefault}{\color[rgb]{0,0,0}$C$}%
}}}}
\put(4426,-3211){\makebox(0,0)[lb]{\smash{{\SetFigFont{17}{20.4}{\familydefault}{\mddefault}{\updefault}{\color[rgb]{0,0,0}$\Omega$}%
}}}}
\put(6451,-4036){\makebox(0,0)[lb]{\smash{{\SetFigFont{17}{20.4}{\familydefault}{\mddefault}{\updefault}{\color[rgb]{0,0,0}$H_2$}%
}}}}
\put(1801,-1231){\makebox(0,0)[lb]{\smash{{\SetFigFont{17}{20.4}{\familydefault}{\mddefault}{\updefault}{\color[rgb]{0,0,0}$\Omega'$}%
}}}}
\put(3601,-211){\makebox(0,0)[lb]{\smash{{\SetFigFont{17}{20.4}{\familydefault}{\mddefault}{\updefault}{\color[rgb]{0,0,0}$H_1$}%
}}}}
\put(3526,-1588){\makebox(0,0)[lb]{\smash{{\SetFigFont{17}{20.4}{\familydefault}{\mddefault}{\updefault}{\color[rgb]{0,0,0}$q_2$}%
}}}}
\put(2176,-5011){\makebox(0,0)[lb]{\smash{{\SetFigFont{17}{20.4}{\familydefault}{\mddefault}{\updefault}{\color[rgb]{0,0,0}$B_{R_\lambda}(p)$}%
}}}}
\end{picture}%
}}
\label{Imma8}\caption{The darker region is $\{ x \in \Omega\, |\, u(x) > h_0\}$.}
\end{figure}

As the length of the segment $\overline{q_{1}q_{2}}$ is bigger or equal
then $2 R_{\lambda}$, and $u(C) \geq h_{0}$, there will exist a first point of
contact between the moved graph $G$ and the graph of $u$ over
$\Omega_{1}$, at the interior or at the boundary of $\Omega$.

Both cases contradict the maximum principle (in the second case because
$\langle \nabla v, \nu \rangle = \alpha$).

We conclude that $d < 2 R_{\lambda}$. The least sentence in the statement
is a classical geometric property which relates the diameter and the
circumradius of a planar figure.

In the case that $\Omega'$ is unbounded there exists a divergent curve $\gamma \subset \Omega'$ and an arc $C \subset \gamma$ whose boundary points are at distance bigger than $2\,R_\lambda$, and we can repeat the previous argument in order to obtain a contradiction. This completes the proof of the result.\proofend

\begin{remark}
If $\Omega$ is an unbounded $f$-extremal domain, and $f$ satisfies property $P_2$, then $\alpha$ cannot be zero. In fact, if $u$ were the solution of (\ref{pr_bis}), then for all $\epsilon$ small enough there would exist an unbounded curve $\Gamma$ in $\Omega$ where $u(\Gamma) > \epsilon$, and one could repeat the argument of the previous proof where $v$ is a solution of (\ref{pr:ball}) such that
$$
\max_{x \in B_{R_\lambda}} v(x) = \epsilon \,,
$$
obtaining a contradiction.
\end{remark}

\section{Boundedness of planar strip ends}

In this section we suppose $\Omega$ to be an unbounded open connected domain of $\mathbb{R}^2$ whose boundary is of class $C^2$, and such that there exists a function $u \in C^2(\overline{\Omega})$ that solves elliptic problem (\ref{pr_bis}), where $f : (0,+\infty) \to \mathbb{R}$ is a Lipschitz function. Moreover we suppose that there exists a constant $R$ such that $\overline{\Omega}$ does not contain any closed ball of radius $R$, i.e., the domain satisfies property $P_1$ (note that this property is satisfied for example if property $P_2$ holds, i.e., when there exists a positive constant $\lambda$ such that $f(t) \geq \lambda\,t$ for all $t>0$, and in this case $R = R_\lambda$). Let $L$ be a straight line of $\mathbb{R}^{2}$, $L^{+}$ and $L^{-}$ the two half-spaces separated by $L$. First we prove a boundedness
property which
will be a key step in the proof of Theorem \ref{T2}. Similar geometric ideas were used
by W. H. Meeks \cite{Meeks} in the context of constant mean curvature surfaces.
For other related boundedness results see the paper of J. A. Aledo, J. M. Espinar and J. A. G\'alvez \cite{al-esp-gal}.

\begin{lemma}\label{bounded}
Let $C$ be a bounded connected component of $\Omega \cap L^+$, and $h(C)$ be the maximum distance of $\pp C$ to $L$. Then 
\[
h(C) \leq 3R. 
\]
\end{lemma}

\textit{Proof.}
Reasoning by contradiction, we will suppose that $h=h(C)$ is greater than $3R$. We can suppose that $L^+$ is the half-space $\{y>0\}$, where $x$ and $y$ denote the coordinates of $\mathbb{R}^2$, and $L=\{y=0\}$.\\
By Proposition \ref{compact} the closure of the curve $\partial C \cap L^+$ is the graph of a function $g(x)$ on a segment of $L$, say $[a,b]$, which is positive at the interior and vanishes at the boundary. Moreover we can assume that the maximum of $g$ is attained at $x=0$, $g(0)=h$, see fig. 4.\\ 
Observe that when one intersects $C$ with the line $\{y=R\}$, the connected components of $C \cap \{y=R\}$ are open intervals whose length is less than $2R$. In fact, if one of such intervals is given by $\{(x,R)\, |\, a' < x < b'\}$ with $b'-a' \geq 2R$, then the rectangle $(a',b')\times(0,R)$ would be contained in $C$ and then $C$ contains a half-ball of radius $R$ and base on the $x$-axis, leading to a contradiction by Corollary \ref{halfAlex}.\\
Let $\tilde C$ be the connected component of $C\cap \{y > R\}$ whose boundary contains the point $(0,h)$. Let $\Gamma$ be the closure of the boundary of $\tilde C$ in $\{y>R\}$ and $p$ and $q$ be the end points of $\Gamma$. Note that $\Gamma$ is a graph over $\{y=R\}$ of height $h_1=h-R$. Moreover $|\overline{pq}|\leq 2R$.
By our hypothesis, $h_1 > 2R$, and then there exists a point $p' \in \Gamma$, other than $q$, maximizing the distance to $p$ and therefore $|\overline{pp'}|> h_1> |\overline{pq}|$, see fig. 4.

\begin{figure}[!ht]
\centering
{\scalebox{.8}{
\begin{picture}(0,0)%
\includegraphics{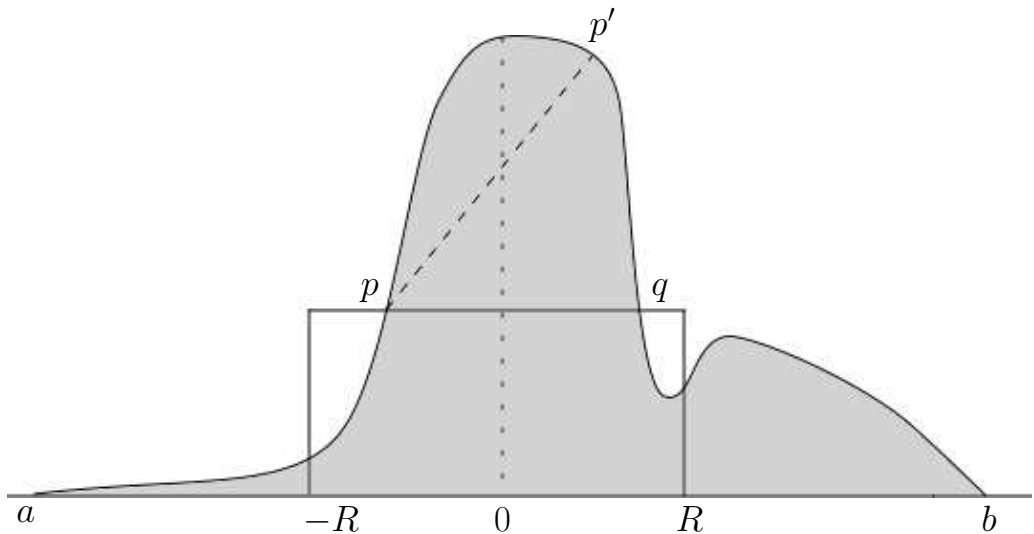}%
\end{picture}%
\setlength{\unitlength}{3947sp}%
\begingroup\makeatletter\ifx\SetFigFont\undefined%
\gdef\SetFigFont#1#2#3#4#5{%
  \reset@font\fontsize{#1}{#2pt}%
  \fontfamily{#3}\fontseries{#4}\fontshape{#5}%
  \selectfont}%
\fi\endgroup%
\begin{picture}(8040,4218)(526,-3748)
\put(8176,-3661){\makebox(0,0)[lb]{\smash{{\SetFigFont{17}{20.4}{\familydefault}{\mddefault}{\updefault}{\color[rgb]{0,0,0}$b$}%
}}}}
\put(5776,-3661){\makebox(0,0)[lb]{\smash{{\SetFigFont{17}{20.4}{\familydefault}{\mddefault}{\updefault}{\color[rgb]{0,0,0}$R$}%
}}}}
\put(4351,-3661){\makebox(0,0)[lb]{\smash{{\SetFigFont{17}{20.4}{\familydefault}{\mddefault}{\updefault}{\color[rgb]{0,0,0}$0$}%
}}}}
\put(2851,-3661){\makebox(0,0)[lb]{\smash{{\SetFigFont{17}{20.4}{\familydefault}{\mddefault}{\updefault}{\color[rgb]{0,0,0}$-R$}%
}}}}
\put(601,-3586){\makebox(0,0)[lb]{\smash{{\SetFigFont{17}{20.4}{\familydefault}{\mddefault}{\updefault}{\color[rgb]{0,0,0}$a$}%
}}}}
\put(3304,-1801){\makebox(0,0)[lb]{\smash{{\SetFigFont{17}{20.4}{\familydefault}{\mddefault}{\updefault}{\color[rgb]{0,0,0}$p$}%
}}}}
\put(5587,-1801){\makebox(0,0)[lb]{\smash{{\SetFigFont{17}{20.4}{\familydefault}{\mddefault}{\updefault}{\color[rgb]{0,0,0}$q$}%
}}}}
\put(5101,239){\makebox(0,0)[lb]{\smash{{\SetFigFont{17}{20.4}{\familydefault}{\mddefault}{\updefault}{\color[rgb]{0,0,0}$p'$}%
}}}}
\end{picture}%
}}
\label{3R}\caption{The bounded component $C$.}
\end{figure}

Denote by $L_*$ the line through $p$ and $p'$ and let $L_*^+$ and $L_*^-$ the half-spaces determined by $L_*$ (i.e., the connected components of $\mathbb{R}^2 \backslash L_*$) such that $q \in L_*^-$. We have that $\tilde C \cap L_*^+$ is a bounded connected component of $\Omega \cap L_*^+$ and by construction it is clear that $L_*$ is orthogonal to the boundary of $\tilde C$ at the point $p'$. By Corollary \ref{cor_orth} we have a contradiction. 
%Apply the Serrin reflection argument starting with a line $L_0 \in L^+$ parallel to $L$ and not intersecting $\tilde C$ (observe that $L_0$ cannot be a vertical line: it should be tilted to the right). Consider the half-spaces $L_t^+$ and $L_t^-$ determined by $L_t$ so that for $t> 0$, we have $p,q\in L_t^+$. Note that in $L_t^-$ we consider only the region $\tilde C$ while in $L_t^+$ one use the whole $\Omega$. That is if $\tilde C_t^- = \tilde C\cap L_t^-$ and $\tilde C_t^*$ is its reflected image with respect to $L_t$, then the reflection argument will continue while $\tilde C_t^*\subset \Omega$. As $L_0$ is normal to $\Gamma$ at $p'$, we get that the symmetry line appears for $t>0$ or at most for $t=0$. Then we conclude that $\Omega$ is symmetric with respect to this $L_t$ and $\Omega$ is a compact region, which is a contradiction.
\proofend

Now we start to study the behavior of $\Omega$ at infinity, where $\Omega$ is supposed to be an $f$-extremal domain satisfying property $P_{1}$ and having finite topology. It is clear that $\Omega$ cannot be the complement of a compact region, and that if $\Omega$ is bounded then it is a ball. The only interesting case of finite topology is then the proper finite topology one, and then we will use the following simple kind of end. 

\begin{definition}\label{endss}
A (planar strip) \textbf{end} of $\Omega$ is an unbounded subdomain $E \subset \overline{\Omega}$, with an homeomorphism $F : [0,1] \times [0,+\infty[ \rightarrow E$ such that :
\begin{enumerate}
	\item $F(0,s) \in \pp \Omega$ for all $s \in[0,+\infty[$,
	\item $F(1,s) \in \pp \Omega$ for all $s \in[0,+\infty[$,
	\item $F(t,s) \in \mathring \Omega$ for all $(t,s) \in ]0,1[ \times [0,+\infty[$.
\end{enumerate}
We will call a \textbf{transversal curve} of the end a curve joining a point of $F\left(\{0\} \times [0,+\infty[\right)$ with a point of $F\left(\{1\} \times [0,+\infty[\right)$ and lying in $E$ (see fig. 5).
\end{definition}

\begin{figure}[!ht]
\centering
{\scalebox{.8}{
\begin{picture}(0,0)%
\includegraphics{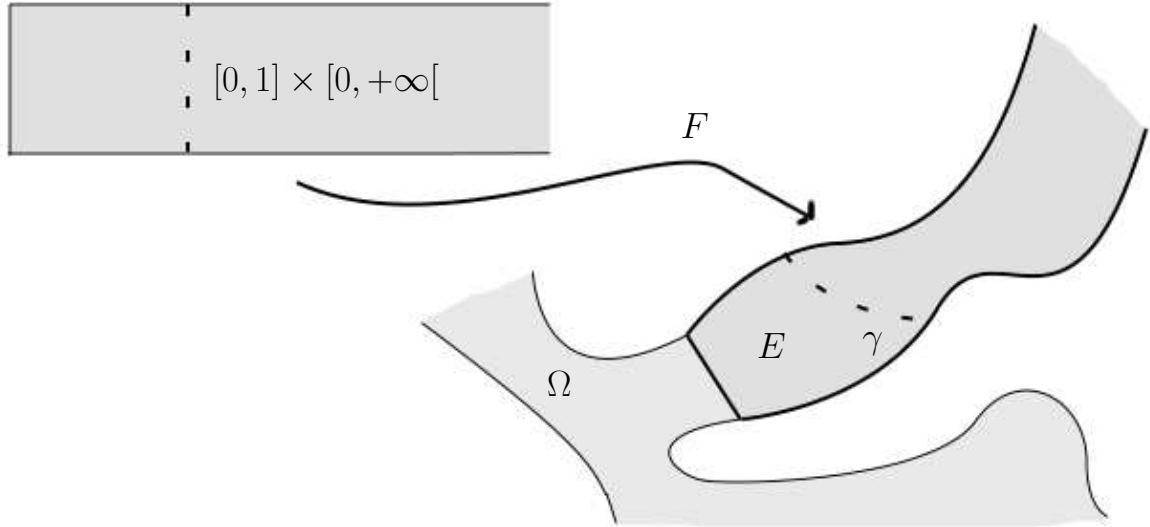}%
\end{picture}%
\setlength{\unitlength}{3947sp}%
\begingroup\makeatletter\ifx\SetFigFont\undefined%
\gdef\SetFigFont#1#2#3#4#5{%
  \reset@font\fontsize{#1}{#2pt}%
  \fontfamily{#3}\fontseries{#4}\fontshape{#5}%
  \selectfont}%
\fi\endgroup%
\begin{picture}(9600,4500)(826,-4336)
\put(6451,-1111){\makebox(0,0)[lb]{\smash{{\SetFigFont{17}{20.4}{\familydefault}{\mddefault}{\updefault}{\color[rgb]{0,0,0}$F$}%
}}}}
\put(7876,-2761){\makebox(0,0)[lb]{\smash{{\SetFigFont{17}{20.4}{\familydefault}{\mddefault}{\updefault}{\color[rgb]{0,0,0}$\gamma$}%
}}}}
\put(7051,-2836){\makebox(0,0)[lb]{\smash{{\SetFigFont{17}{20.4}{\familydefault}{\mddefault}{\updefault}{\color[rgb]{0,0,0}$E$}%
}}}}
\put(5401,-3136){\makebox(0,0)[lb]{\smash{{\SetFigFont{17}{20.4}{\familydefault}{\mddefault}{\updefault}{\color[rgb]{0,0,0}$\Omega$}%
}}}}
\put(2776,-736){\makebox(0,0)[lb]{\smash{{\SetFigFont{17}{20.4}{\familydefault}{\mddefault}{\updefault}{\color[rgb]{0,0,0}$[0,1]\times[0,+\infty[$}%
}}}}
\end{picture}%
}}
\label{end}\caption{An end $E$ and a transversal curve $\gamma$}
\end{figure}

Let $E$ be a (planar strip) end of $\Omega$ and let $L$ be a straight line of $\mathbb{R}^2$ intersecting $E$. The first property of the ends of our domain is the following:

\begin{lemma}\label{nonbounded}
Let $E$ be an end of $\Omega$ and suppose that $L \cap E$ contains an unbounded connected component. Then any straight line $L'$ parallel to $L$ and sufficiently far from $L$ intersects $E$ in only bounded connected components.
\end{lemma}

\textit{Proof.} After a rigid motion, we can suppose that $L$ is the $x$-axis of $\mathbb{R}^2$ and the unbounded connected component of $L \cap E$ is
\[
\{(x,0) \in \mathbb{R}^2\, |\, x \in [0,+\infty[ \,\, \, \}
\] 
Now take a straight line $L'$, parallel to $L$ and at a distance from $L$ bigger than $R$. Of course $L'$ is given by the equation $y=k$ with $|k| > R$. If $L' \cap E$ contains an unbounded connected component, then there exists a constant $\rho$ such that the unbounded connected component $C$ of $L \cap E$ is either
\[
\{(x,k) \in \mathbb{R}^2\, |\, x \in [\rho,+\infty[\,\, \,\}\qquad \textnormal{or}\qquad \{(x,k) \in \mathbb{R}^2\, |\, x \in ]-\infty, \rho]\,\,\}
\] 
Moreover, there exists a regular curve $\gamma \in \mathbb{R}^2$ joining $(0,0)$ to $(\rho,k)$ and lying in $E$. We have that 
\[
\sigma = \{(x,0) \in \mathbb{R}^2\, |\, x \in [0,+\infty[\,\, \,\} \,\, \, \cup \,\, \, \gamma \,\, \, \cup \,\, \, C
\]
separates $\mathbb{R}^2$ in two connected components, one of which is contained in $E$, and then also in $\Omega$. Lemma \ref{balls_R} leads to a contradiction because both the components of $\mathbb{R}^2 \backslash \sigma$ contain balls of radius $R$. Hence $L' \cap E$ does not contain any unbounded connected component and the lemma follows at once.
\proofend

The main result of this section is Statement (T1) of Theorem \ref{T2}:

\begin{proposition}\label{bounded_end}
Let $E$ be a (planar strip) end of $\Omega$. Then $E$ stays at bounded distance from a half-line.
\end{proposition}

\textit{Proof.}  Let $F$ be the homeomorphism associated to $E$ by definition \ref{endss}, between
$[0,1]\times [0,\infty[$ and $E$, and
$\beta= F([0,1]\times \{0\})$ the initial transversal
curve of the end. Let $B=B_r(0)$ be a ball of radius
$r$ centered at the origin of $\mathbb{R}^2$ containing
$\beta$ and let $p_1, p_2, p_3, \ldots$ be a divergent sequence
of points in $E$ such that the sequence of normalized
vectors $q_i = p_i/|p_i|$ converges to a unit vector $q$.
After a possible rotation of $E$ we can assume $q = (1, 0)$.\\
We show now that $E$ stays at bounded distance from the
$x$-axis. Otherwise, assume that $E$ intersects every
horizontal line in $y>0$. Choose $\alpha > r+1$
such that $l_\alpha =\{y=\alpha\}$ meets $\partial E$
transversally. By Proposition \ref{compact} and Lemma \ref{bounded}, the region
$E \cap \{y > \alpha\}$ has an unbounded connected
component $C$ and by Lemma \ref{nonbounded} the intersection of $E$
with the line $l_\alpha$ does not contain any
unbounded connected component. Therefore there exists
a transversal curve $\gamma$ of $E$ contained in $C$.\\
It follows that the transversal curves $\beta$ and $\gamma$
lie below and above $l_\alpha$, respectively.

\begin{figure}[!ht]
\centering
{\scalebox{.8}{%\input{image6.pstex_t}
\begin{picture}(0,0)%
\includegraphics{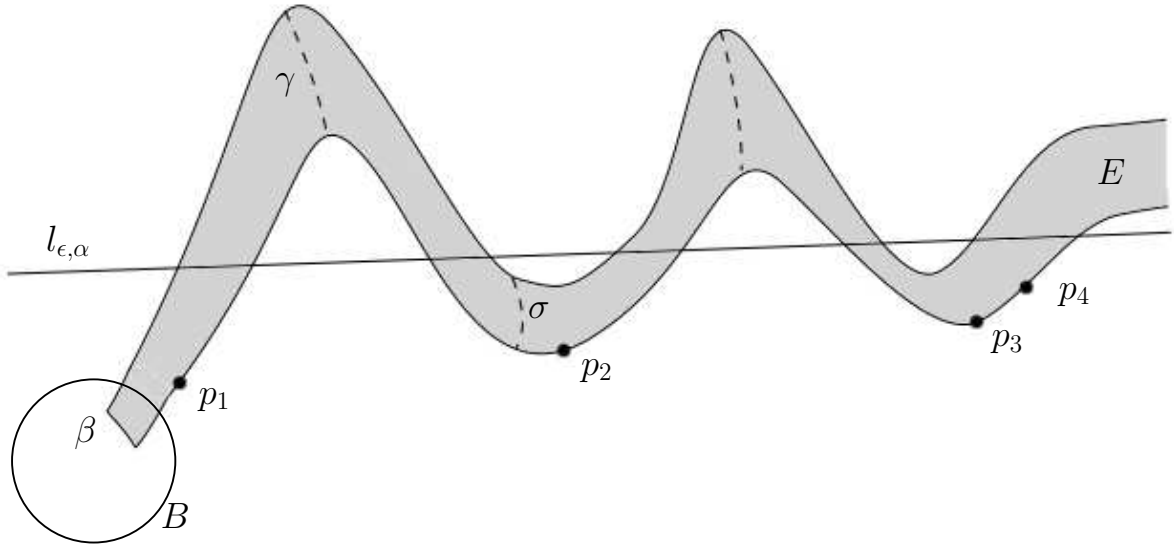}%
\end{picture}%
\setlength{\unitlength}{3947sp}%
\begingroup\makeatletter\ifx\SetFigFont\undefined%
\gdef\SetFigFont#1#2#3#4#5{%
  \reset@font\fontsize{#1}{#2pt}%
  \fontfamily{#3}\fontseries{#4}\fontshape{#5}%
  \selectfont}%
\fi\endgroup%
\begin{picture}(9600,4500)(676,-4261)
\put(3001,-661){\makebox(0,0)[lb]{\smash{{\SetFigFont{17}{20.4}{\familydefault}{\mddefault}{\updefault}{\color[rgb]{0,0,0}$\gamma$}%
}}}}
\put(9151,-2311){\makebox(0,0)[lb]{\smash{{\SetFigFont{17}{20.4}{\familydefault}{\mddefault}{\updefault}{\color[rgb]{0,0,0}$p_4$}%
}}}}
\put(9451,-1411){\makebox(0,0)[lb]{\smash{{\SetFigFont{17}{20.4}{\familydefault}{\mddefault}{\updefault}{\color[rgb]{0,0,0}$E$}%
}}}}
\put(4987,-2461){\makebox(0,0)[lb]{\smash{{\SetFigFont{17}{20.4}{\familydefault}{\mddefault}{\updefault}{\color[rgb]{0,0,0}$\sigma$}%
}}}}
\put(2401,-3136){\makebox(0,0)[lb]{\smash{{\SetFigFont{17}{20.4}{\familydefault}{\mddefault}{\updefault}{\color[rgb]{0,0,0}$p_1$}%
}}}}
\put(5401,-2911){\makebox(0,0)[lb]{\smash{{\SetFigFont{17}{20.4}{\familydefault}{\mddefault}{\updefault}{\color[rgb]{0,0,0}$p_2$}%
}}}}
\put(1201,-1936){\makebox(0,0)[lb]{\smash{{\SetFigFont{17}{20.4}{\familydefault}{\mddefault}{\updefault}{\color[rgb]{0,0,0}$l_{\epsilon,\alpha}$}%
}}}}
\put(8626,-2686){\makebox(0,0)[lb]{\smash{{\SetFigFont{17}{20.4}{\familydefault}{\mddefault}{\updefault}{\color[rgb]{0,0,0}$p_3$}%
}}}}
\put(1426,-3436){\makebox(0,0)[lb]{\smash{{\SetFigFont{17}{20.4}{\familydefault}{\mddefault}{\updefault}{\color[rgb]{0,0,0}$\beta$}%
}}}}
\put(2101,-4111){\makebox(0,0)[lb]{\smash{{\SetFigFont{17}{20.4}{\familydefault}{\mddefault}{\updefault}{\color[rgb]{0,0,0}$B$}%
}}}}
\end{picture}%
}}
\label{fig5}\caption{The end $E$.}
\end{figure}

Therefore, when $\varepsilon >0$
is small enough,  the same holds for the line
$l_{\alpha,\varepsilon}= \{ y=\varepsilon \, x + \alpha\}$,
that is $\beta\subset \{ y < \varepsilon \, x + \alpha\}$
and $\gamma \subset \{ y > \varepsilon \, x + \alpha\}$, see fig. 6.\\
Now we apply again the same argument to the subend $E^*$
of $E$ whose initial transversal arc is $\gamma$, i.e., $E^*$ is the closure of the unbounded component of $E \backslash \gamma$. As the
arc $\gamma$ lies above $l_{\alpha,\varepsilon}$,
almost all points $p_i$ belong to $\{y < \varepsilon \, x + \alpha\}$, and the distance between $p_i$ and
$l_{\alpha,\varepsilon}$ diverges to infinity. Reasoning as above we find that there exists a
transverse curve $\sigma$ of $E^*$ (and so of $E$ also)
contained in $y < \varepsilon \, x + \alpha$ (see fig. 6).\\
The existence of $\gamma$ and $\sigma$ leads to a contradiction. In fact, the component of $E\cap\{y > \varepsilon \, x + \alpha\}$
containing the transversal arc $\gamma$ must be bounded, which
contradicts Proposition \ref{compact} and then $E$ stays at bounded distance from the $x$-axis. \\
In order to prove that $E$ is at bounded distance from the half-line $\{y=0,\, x>0\}$, let
$\Gamma = \partial E \cap \partial \Omega$ and take $b>0$.
Assume that $B\cap \{x< -b\} = \emptyset$ and the line $\{x = -b\}$ intersects $\Gamma$
transversally. Hence $E \cap \{x = -b\}$ consists of a finite union of proper
embedded arcs whose extremes belongs to $\Gamma$. The existence of the divergent sequence $p_{i} \in E \cap \{x > 0, |y| < k\}$, where $k$ is the maximum distance of the end $E$ to the $x$-axis, implies that $E \cap \{x < -b\}$ has only bounded components and using Lemma \ref{bounded} we conclude that $E$ is contained is the half-strip $\{x> -(b+3R), \, |y| < k\}$. Hence $E$ is at bounded distance from a half-line and the proposition follows.\\
\proofend

Now we are able to prove Statements (T2) and (T3) of Theorem \ref{T2}. 

\begin{proposition} \label{main} The following properties hold:
\begin{enumerate}
	\item $\Omega$ cannot have only one (planar strip) end. Moreover $\Omega$ cannot stay in a half-strip.
	\item If $\Omega$ has exactly two (planar strip) ends, then there exists a line $L$ such that $\Omega$ stays at bounded distance from $L$, and the two ends are on opposite sides with respect to any line orthogonal to $L$.
\end{enumerate}
\end{proposition}

\textit{Proof.} The proof of the two statements follows from Propositions \ref{bounded_end}, \ref{compact} and Lemma \ref{bounded}.
\begin{enumerate}
	\item Let us suppose that $\Omega$ is contained in a half-strip. We can suppose that 
\[
\Omega \subseteq \{(x,y) \in \mathbb{R}^2 \,\,:\,\, -A < x < A, y > 0\}
\]
for some positive constant $A$. For any $k>0$, each connected component of $\Omega\cap\{y < k\}$ is bounded. Choosing $k$ large enough one leads to a contradiction by Lemma \ref{bounded}. For a general statement of this kind see Proposition \ref{wedge} below.\\
Let us suppose now that $\Omega$ has only one end $E$. This means that $\Omega \backslash E$ is bounded. By Proposition \ref{bounded_end}, $\Omega$ lies in a half-strip, contradiction.\\
	\item Let us suppose that $\Omega$ has exactly two ends $E_1$ and $E_2$. By Proposition \ref{bounded_end} $E_1$ is at bounded distance from a line $L_1$ and $E_2$ is at bounded distance from an other line $L_2$. If $L_1$ and $L_2$ are parallel, then it is clear that $\Omega$ is at bounded distance from both $L_1$ and $L_2$, because $\Omega \backslash (E_1 \cup E_2)$ is bounded. If $L_1$ and $L_2$ are not parallel, then there exists a straight line $l$ such that $l$ is tangent to $\partial\Omega$ and $\Omega$ is contained in only one of the two connected components of $\mathbb{R}^2 \backslash l$. After a rigid motion we can suppose that $l$ is the $x$-axis and $\Omega$ stays in the upper half-plane. Proposition \ref{bounded_end} implies that each connected component of $\Omega \cap \{y<k\}$, $k>0$, is bounded, and by construction the distance of $\pp \Omega \cap \{y<k\}$ is equal to $k$. Choosing $k$ big enough one leads to a contradiction by Lemma \ref{bounded}. This proves that $L_1$ and $L_2$ must be parallel, and then $\Omega$ is
at bounded distance from a line $L$.
If the two ends are on the same side with respect to a line orthogonal to
$L$, then $\Omega$ is contained in a half-strip, contradiction. Then the two ends are on opposite sides with respect to any line orthogonal to $L$. 
	\end{enumerate}
\proofend

\section{Boundedness of the domains}

In this section we prove some properties of planar domains $\Omega$ where problem (\ref{pr_bis}) can be solved, without any extra assumption on the function $f$, that is only supposed to have Lipschitz regularity. From now to the end of the section $\Omega$ will be a planar $C^2$-domain where problem (\ref{pr_bis}) can be solved.

\medskip

In the case that $\Omega$ is unbounded, and bounded by a unique proper arc $\Gamma$, let
$\gamma(t)=(x(t),y(t))$, $t\in\mathbb{R}$, be an arc-length parametrization of $\Gamma$, 
$\gamma'=(x',y')$ the unit tangent vector and $n=-\nu=(-y',x')$ the inward pointing normal vector along 
$\Gamma$. A basic property of the domain $\Omega$ is the following:

\begin{lemma}\label{halfline}
If $\Omega$ is unbounded, for any point $p\in \Gamma$, the normal inward half-line 
\[
L^+(p)=\{p+t \, n(p)\, | \, t\geq 0 \}
\]
lies in $\{p\}\cup \Omega$ (see fig. 7).
\end{lemma}

\begin{figure}[!ht]
\centering
{\scalebox{.8}{
\begin{picture}(0,0)%
\includegraphics{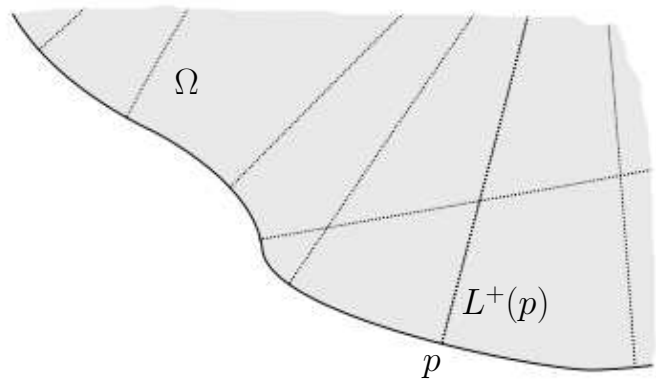}%
\end{picture}%
\setlength{\unitlength}{3947sp}%
\begingroup\makeatletter\ifx\SetFigFont\undefined%
\gdef\SetFigFont#1#2#3#4#5{%
  \reset@font\fontsize{#1}{#2pt}%
  \fontfamily{#3}\fontseries{#4}\fontshape{#5}%
  \selectfont}%
\fi\endgroup%
\begin{picture}(5745,3240)(226,-2626)
\put(2026,-286){\makebox(0,0)[lb]{\smash{{\SetFigFont{17}{20.4}{\familydefault}{\mddefault}{\updefault}{\color[rgb]{0,0,0}$\Omega$}%
}}}}
\put(4276,-2011){\makebox(0,0)[lb]{\smash{{\SetFigFont{17}{20.4}{\familydefault}{\mddefault}{\updefault}{\color[rgb]{0,0,0}$L^{+}(p)$}%
}}}}
\put(3976,-2461){\makebox(0,0)[lb]{\smash{{\SetFigFont{17}{20.4}{\familydefault}{\mddefault}{\updefault}{\color[rgb]{0,0,0}$p$}%
}}}}
\end{picture}%
}}
\label{fig7}\caption{The normal inward half-line $L^+(p)$, being $p$ a point of $\partial \Omega$.}
\end{figure}

{\it Proof.} In fact, this holds for small $t>0$ and if $L^+(p)$ meets $\Gamma$ in a second point $p'$, for the first time, then there exists a bounded region $C\subset\Omega$ bounded by the segment 
$\overline{pp'}$ and the arc in $\Gamma$ joining $p$ and $p'$. As both arcs cut orthogonally at $p$, 
we have a contradiction by using the moving line argument as in the proof of Lemma \ref{bounded} (see fig. 4 and invert $p$ and $p'$).
\proofend

We are now able to prove Statements (T6) and (T7)  of Theorem \ref{T4}.

\begin{proposition} \label{wedge}
If $\Omega$ is contained in a wedge of angle less than $\pi$ 
(no restriction about the topology of the domain), then $\Omega$ is either a ball or a uniformly Lipschitz epigraph. 
If the angle of the wedge is less than $\pi/2$, then the domain is a ball.
\end{proposition}

{\it Proof.} If we choose the coordinates so that the wedge is contained in the upper half-plane and is symmetric with respect to $y$-axis, then $\Omega \cap \{y<a\}$ is bounded for any $a>0$ and using the Serrin's reflection argument with horizontal lines, as in the proof of Proposition \ref{main}, we conclude that either $\Omega$ is bounded (and then a ball, by the theorem of J. Serrin) or $\Gamma=\partial \Omega$ is a proper arc whose projection over the $x$-axis is one-to-one. In this last case, the previous lemma asserts that $L^+(p)$ lies in the wedge for all $p \in \Gamma$ and it follows that 
$x'>\varepsilon$ for some positive $\varepsilon$, where $(x(t),y(t))$, $t\in\mathbb{R}$, is an arc-length parametrization of $\Gamma$. Hence $\Omega$ is a uniformly Lipschitz epigraph. \\
If the angle of the wedge is smaller than $\pi/2$, then by changing the Euclidean coordinates, we can assume that $\Omega \subset \{ 0< y < bx \}$  for some $b>0$. Lemma \ref{halfline} implies that, for any point $p\in\Gamma$, the slope of the normal half-line $L^+(p)$ lies between the ones
of $\{y=0\}$ and $\{y=bx\}$ and therefore, the unit tangent vector $\gamma'(t)$, $t\in\mathbb{R}$, points down more than the vector $(b,-1)$ which contradicts the fact that $\Omega \subset \{y>0\}$.  
\proofend

\begin{remark}
Let $\Gamma$ be the graph of a $C^2$ function
$\varphi:\mathbb{R}\longrightarrow\mathbb{R}$ and $\Omega$
the epigraph $\{y>\varphi(x)\}$. Let's see some cases where
$\Omega\subset\mathbb{R}^2$ satisfies the conclusion of
Lemma \ref{halfline}, i.e. for any $p\in \Gamma$,
\begin{equation}
\label{property}
L^+(p)\subset  \Omega \cup \{p\}.
\end{equation}

\vspace{.3cm}

 The tangent vector at a point $p=(x,\varphi(x))$ is
$(1,\varphi'(x))$ and the inner normal half-line is
$
L^+(p)=\{(x,\varphi(x))+a(-\varphi'(x),1) \, / \, a\geq 0\}.
$

\vspace{.3cm}

1) If $\varphi'\geq 0$, then $\varphi$ is increasing and
the inner normal half-line is either vertical or tilted to
the left. Then (\ref{property}) follows.

\vspace{.2cm}

2) If $\varphi''\leq 0$, the $\mathbb{R}^2-\Omega$ is
convex and so $\Omega$ satisfies (\ref{property}). On
the contrary, if $\Omega$ is convex, then (\ref{property})
does not hold in general. However it can be verified
directly in some cases like for the domain bounded by the
equilateral hyperbola $\Omega=\{y>\sqrt{1+x^2}\}$. In
particular, by using our argument we cannot improve
Proposition \ref{wedge} to include the case $\theta=\pi/2$ in the second statement.

\vspace{.2cm}

3) If $|\varphi'| \leq 1$, then the epigraph $\Omega$ satisfies
(\ref{property}). Otherwise, we can suppose there is
$x_1<x_2$ such that the inner normal half-line $L^+(p)$,
with $p=(x_1,\varphi(x_1))$, meets $\Gamma$ at the point
$q=(x_2,\varphi(x_2))$. Therefore, the slope of the
$L^+(p)$ is positive and we have that
\[
1\leq \frac{-1}{\varphi'(x_1)}=\frac{\varphi(x_2)-
\varphi(x_1)}{x_2-x_1}=\frac{1}{x_2-x_1}\int_{x_1}^{x_2}
\varphi' \, dx\leq 1.
\]
It follows that both inequalities are in fact equalities:
$\varphi'(x_1)=-1$ and $\varphi'(x)=1$, $x_1\leq x \leq x_2$
and this contradiction proves the claim.
For instance, property (\ref{property}) is satisfied for
the functions
\[
\varphi(x)= \sin(x), \hspace{1cm}
\varphi(x)=\frac{1}{4}\log(1+x^2)\, \sin(\log(1+x^2)).
\]
The first function is periodic and the second one is oscillating
and gives a domain $\Omega$ which neither contains a half-plane nor
is contained in a half-plane; compare with Proposition \ref{cor2BCN}.
\end{remark}

Now we prove Statement (T8)  of Theorem \ref{T4}. The basic idea of the proof is the tilted moving plane argument, used in \cite{KKS} for surface theory. 

%\begin{proposition}
%If $\Omega$ is contained in a half-plane and its boundary consists of a unique open curve $\Gamma$, then $\Omega$ is an weak epigraph.
%\end{proposition}

%{\it Proof.} Assume that $\Omega\subset \{y>0\}$. First observe that Lemma \ref{halfline} implies that, for any point $p \in \Gamma$, the normal half-line $L^+(p)$ lies in $\{y>0\}$. 
%As a consequence, $x'\geq 0$ for any real $t$. Otherwise, for some $p\in \Gamma$, the normal inward vector $n$ at $p$ points down and the normal half-line $L^+(p)$ would enter in the half-plane below the $x$-axis. Then, the projection of $\Gamma$ over the $x$-axis is parametrized by the function $x(t)$ which is monotonically increasing. If the image of $x$ is the whole $x$-axis, then the statement follows. If the supremum of the projection of $\Gamma$ over the $x$-axis is a real number $a$, then the  domain is contained in the wedge $\{ x<a, y>0\}$ and it follows from Proposition \ref{wedge} that the domain is an epigraph. The same argument applies if the infimum of the projection of $\Gamma$ is finite. This completes the proofs of the proposition. \proofend

\begin{proposition}\label{halfplane}

Let $\Omega\subset \mathbb{R}^2$  be an $f$-extremal domain
contained in a half-plane. Then $\Omega$ is either a ball or (after
a rigid motion) there exist a $C^2$ positive function
$\varphi:\mathbb{R}\longrightarrow ]0,\infty[$ such that either

\begin{enumerate}
\item[i.]the domain $\Omega$ is an epigraph $\{y>\varphi (x)\}$, or
\item[ii.] $\varphi$ is bounded and $\Omega$ is the symmetric domain
$\{|y|<\varphi(x)\}$.
\end{enumerate}

\end{proposition}

{\it Proof.} If $\Omega$ is bounded, then it is a ball by Serrin's theorem. Then, assume that $\Omega$ is an unbounded subset of the half-plane $\{y>0\}$. After a suitable translation, we can assume
that $\partial\Omega$ intersects the $y$-axis transversally.\\
Let us suppose that the intersection of $\partial\Omega$ with the $y$-axis is done by more than one point. Consider $\Omega_1=\Omega\cap \{x>0\}$
and $\Omega_2=\Omega\cap \{x<0\}$. Note that $\Omega_1$ and
$\Omega_2$ are nonempty open sets. If either $\Omega_1$ or $\Omega_2$ is contained in a vertical slab, then it follows from Proposition \ref{wedge} that $\Omega$ is an epigraph with
respect to one of the diagonals of the plane. So, henceforth we will suppose that $\Omega_1$ and $\Omega_2$ are unbounded
and that the orthogonal projection of $\partial \Omega$ over $\{y=0\}$ is onto. It is clear that $\Omega\cap \{x=0\}$
is a discrete union of open intervals in the $y$-axis,
the lowest of these intervals being bounded. Denote by
$p$ the lowest point in the boundary of this interval.\\
Given a straight line $T$, for any $x\in\mathbb{R}^2$ and any
subset $X\subset\mathbb{R}^2$ let $x'$ be the reflection of $x$ about $T$
and $X'$ be the reflected image of $X$ about $T$. Fix $\varepsilon>0$ and consider the two pencils of parallel straightlines
\[
T_a =\{ y= a\} \qquad \textnormal{and} \qquad T_{\varepsilon,a}=
\{y=-\varepsilon \, x+ a\}
\]
for $a \in \mathbb{R}$. Now we use the moving line argument. Let $T=T_{\varepsilon,a}$
be an element of the second pencil. For $a=0$ the line $T$
does not intersect $\Omega_1$. We suppose this line to be
continuously moved parallel to itself, by increasing $a$,
until it pass through $p$. From that moment onward, at each
stage of the motion the resulting line $T$ will cut off from
$\Omega_1$ a bounded cap $\Sigma(T)$ defined as follows.
As the part of $\Omega_1$ below $T$ is bounded, it follows
from Proposition \ref{compact} that the reflected image with respect
to $T$ of the connected components of
$\Omega_1\cap\{y<-\varepsilon \, x +a\}$ are contained in
$\Omega$, except possibly for the component whose boundary
contains $p$. Let's denote this component by $\Sigma(T)$.
The portions of the boundary of $\Sigma(T)$ contained in
$T$, $x=0$ and $\partial\Omega$ will be denoted by $I$,
$J$ and $K$, respectively. Note that $p\in J\cap K$.\\
Let  $\Sigma'(T)$, $K'$ and $J'$ be respectively the symmetric image of  $\Sigma(T)$, $K$ and $J$ about $T$. Define on the closure of $\Sigma'(T)$
the function $u'_T$ given by $u'_T (x) = u(x')$.
At the beginning $\Sigma'(T)$ is contained in $\Omega$
and $u'_T\leq u$ and we continue the process while this
%while $\Sigma'(T)\subset \Omega$
occurs. As the $y$-axis cuts transversally $\partial \Omega$
in at least two points, we will meet a first value
$a=a(\varepsilon)>0$  for which one of the following
events holds (see fig. 8):

\vspace{.2cm}

$(1)$  at an interior point, the reflected arc $K'$ touches
the boundary of $\Omega$,

\vspace{.2cm}

$(2)$ $K$ meets $T$ orthogonally,

\vspace{.2cm}

$(3)$ at a point of $\Sigma'(T)\cup I$, the graph
of the resulting function $u'_T$ is tangent to the graph
of the function $u$,

\vspace{.2cm}

(4) $p'$, the reflection of $p$ about $T$,
belongs to $\partial\Omega$,

\vspace{.2cm}

(5) when restricted to the segment $J'$, the graph of
the resulting function  $u'_T$ is tangent at some
interior point to the graph of the function $u$.

\begin{figure}[!ht]
\centering
{\scalebox{.8}{%\input{image8.pstex_t}
\begin{picture}(0,0)%
\includegraphics{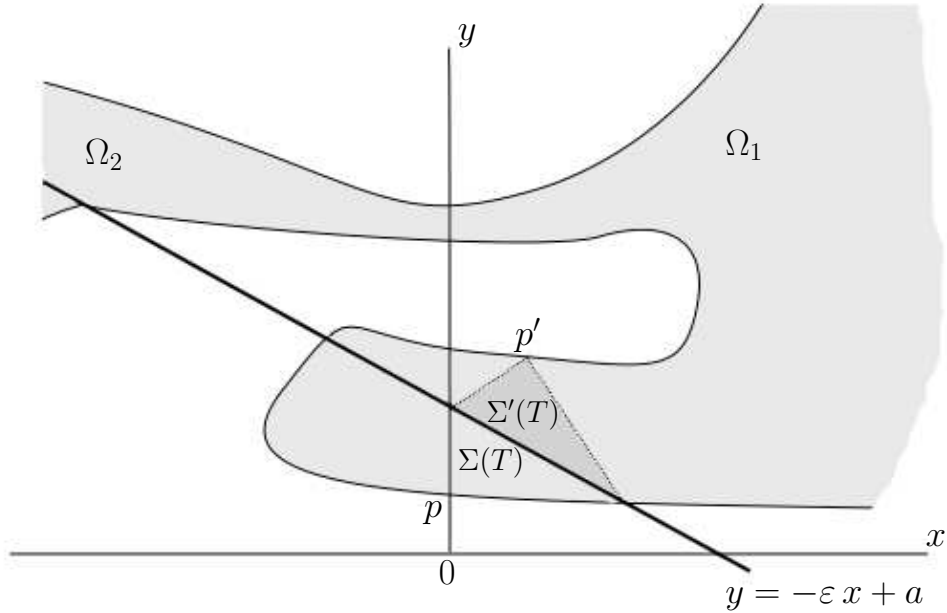}%
\end{picture}%
\setlength{\unitlength}{3947sp}%
\begingroup\makeatletter\ifx\SetFigFont\undefined%
\gdef\SetFigFont#1#2#3#4#5{%
  \reset@font\fontsize{#1}{#2pt}%
  \fontfamily{#3}\fontseries{#4}\fontshape{#5}%
  \selectfont}%
\fi\endgroup%
\begin{picture}(7680,5070)(226,-4456)
\put(6001,-811){\makebox(0,0)[lb]{\smash{{\SetFigFont{17}{20.4}{\familydefault}{\mddefault}{\updefault}{\color[rgb]{0,0,0}$\Omega_1$}%
}}}}
\put(976,-886){\makebox(0,0)[lb]{\smash{{\SetFigFont{17}{20.4}{\familydefault}{\mddefault}{\updefault}{\color[rgb]{0,0,0}$\Omega_2$}%
}}}}
\put(3901, 89){\makebox(0,0)[lb]{\smash{{\SetFigFont{17}{20.4}{\familydefault}{\mddefault}{\updefault}{\color[rgb]{0,0,0}$y$}%
}}}}
\put(4351,-2311){\makebox(0,0)[lb]{\smash{{\SetFigFont{17}{20.4}{\familydefault}{\mddefault}{\updefault}{\color[rgb]{0,0,0}$p'$}%
}}}}
\put(7576,-3886){\makebox(0,0)[lb]{\smash{{\SetFigFont{17}{20.4}{\familydefault}{\mddefault}{\updefault}{\color[rgb]{0,0,0}$x$}%
}}}}
\put(3901,-3286){\makebox(0,0)[lb]{\smash{{\SetFigFont{14}{16.8}{\familydefault}{\mddefault}{\updefault}{\color[rgb]{0,0,0}$\Sigma(T)$}%
}}}}
\put(4126,-2911){\makebox(0,0)[lb]{\smash{{\SetFigFont{14}{16.8}{\familydefault}{\mddefault}{\updefault}{\color[rgb]{0,0,0}$\Sigma'(T)$}%
}}}}
\put(3751,-4186){\makebox(0,0)[lb]{\smash{{\SetFigFont{17}{20.4}{\familydefault}{\mddefault}{\updefault}{\color[rgb]{0,0,0}$0$}%
}}}}
\put(3637,-3661){\makebox(0,0)[lb]{\smash{{\SetFigFont{17}{20.4}{\familydefault}{\mddefault}{\updefault}{\color[rgb]{0,0,0}$p$}%
}}}}
\put(6001,-4336){\makebox(0,0)[lb]{\smash{{\SetFigFont{17}{20.4}{\familydefault}{\mddefault}{\updefault}{\color[rgb]{0,0,0}$y=-\varepsilon\,x + a$}%
}}}}
\end{picture}%
}}
\label{fig8}\caption{Tilted moving plane method.}
\end{figure}

By the Serrin's reflection method, we deduce that each one
of the first three options implies that 
$K'\subset\partial \Omega$.
Therefore both events $(4)$ and $(5)$ are also true.
We conclude that in fact the process can be carried on until either event (4) or event (5) occurs for a first value
$a=a(\varepsilon)>0$. \\ 
Now take a sequence of $\varepsilon_i>0$ going to zero, and repeat all the reasoning with $\varepsilon = \varepsilon_i$. The sequence $a(\varepsilon_i)$ is bounded and then, if $a$ is the limit of $a(\varepsilon_i)$, the argument and its conclusion hold also for the limit
horizontal line, leading to the following result:
there exists a horizontal line $T = T_a$, with
$a>0$, such that the reflected image of
$\Omega_1 \cap \{y < a\}$ lies in $\Omega$,
$u'_T\leq u$ and one of the two events $(4)$ or $(5)$
above occurs. Moreover, as $\varepsilon=0$, $J$ is an
interval contained  in the closure of the lowest
interval of $\Omega\cap\{x=0\}$ and the value of $a$
depends only on $J$ and on the behavior of $u$
restricted to $\Omega\cap\{x=0\}$.\\
Now repeat all the process for
$\Omega_2=\Omega\cap\{y < 0\}$ instead of
$\Omega_1$, with lines of positive slope defined
by $T^{*}_{\varepsilon,a}=\{y=\varepsilon x +a\}$.
We obtain the existence of a horizontal line
$T^* = \{y = a^*\}$, such that the reflected image of
$\Omega_2 \cap \{ y < a^*\}$ stays within $\Omega$,
$u'_T\leq u$ and one of the two events $(4)$ or
$(5)$ occurs.\\
As $a$ and $a^*$ depends only on the behavior
of the solution $u$ along $x=0$, it follows that $a=a^*$
and the line $T=T^*$ satisfies that the reflected image of
$\Omega\cap \{y<a\}$ with respect to $T$ is contained
in $\Omega$, $u'_T\leq u$ and one of the assertions
$(1)$, $(2)$ or $(3)$ holds (at some point of the $y$-axis).
From the Serrin's reflection argument we obtain that $\Omega$
is symmetric with respect to $T$. After a suitable rigid motion, item
$ii)$ in the statement of the proposition follows for the domain $\Omega$ from the fact that $\Omega_1$ and $\Omega_2$ are both unbounded.\\
Now let us consider the case where any vertical line which meets transversally the boundary of $\Omega$
meets $\partial\Omega$ just in a point. Then the boundary
of $\Omega$ consists of a unique proper arc $\Gamma$
which projects monotonically and surjectively
onto the $x$-axis.
If $\Gamma$ is given as the graph of a function, then
$\Omega$ is an epigraph.
If the arc $\Gamma$ is tangent to a vertical line at some point
$q$ of the horizontal line $T=\{y=b\}$, $b>0$, then we repeat
the reflection argument of the beginning of the proof with
straight lines $T = T_{\epsilon,a}$ and $T^*=T^*_{\epsilon,a}$, $a\leq b$, and we conclude
that the domain is symmetric with respect to a line $T=\{y=a\}$, with $a\leq b$, and this is not possible by the assumptions on $\Gamma$.
This contradiction completes the proof of the proposition.
\proofend

Statement (T9) of Theorem \ref{T4}.  is a consequence of the previous proposition. Its proof follows from Corollary \ref{cor_BCN} and the following:

\begin{proposition}\label{cor2BCN}
Let $\Omega$
be an $f$-extremal unbounded domain of $\mathbb{R}^2$ bounded by a unique proper arc. Then
either $\Omega$ is an almost epigraph or it contains a half-plane.
\end{proposition}

{\it Proof.} Let $\Gamma$ be the boundary of $\Omega$ and
$\gamma(t)=(x(t),y(t))$, $t\in\mathbb{R}$, an arc-length parametrization of $\Gamma$, 
$\gamma'=(x',y')$ the unit tangent vector and $n=-\nu=(-y',x')$ the inward pointing normal vector along  $\Gamma$. If there are two points $p,q\in\Gamma$ such that $n(p)=-n(q)$,
then from Lemma \ref{halfline} we have that $\Omega$ contains two parallel
half-lines with opposite orientation and it follows that $\Omega$
contains a half-plane. Otherwise the normal image is contained in
a half-circle and, after a rigid motion we can assume that
$x'\geq0$. If the image $I$ of $t \to x(t)$ coincides with $\mathbb{R}$ then $\Omega$ is an almost epigraph. If $I$ is not the whole $x$-axis, we can suppose for example that it is bounded above by a constant $a$, and therefore
$\Gamma$ and $x=a$ are disjoint and then $\Omega$ is contained in a half-plane. Using Proposition \ref{halfplane}, we get that $\Omega$ is an epigraph.\proofend

The first part of the proof of Proposition \ref{halfplane} can be adapted to $f$-extremal domains $\Omega$ in dimension $n \geq 3$ contained in a solid cylinder in order to obtain Theorem \ref{T0}.

\medskip

\textit{Proof of Theorem \ref{T0}.} Suppose that the unbounded domain $\Omega \subset \mathbb{R}^n$ is contained in the solid cylinder $\mathbb{R} \times B$, where $B$ is a $(n-1)$-dimensional ball.
If $\Omega$ is contained in a half-cylinder, assume $\Omega \subset \{x>0\} \times B$, then by cutting the figure with hyperplanes normal to the $x$-axis, and using
Proposition \ref{compact}, we conclude that $\partial\Omega$ is the graph of a
function $f:D\to\mathbb{R}$, with $D \subset B$, such that the limit value of
$f$ at $\partial D$ is $+\infty$. So it follows that a suitably chosen
half-plane $H^- = \{y_1 \geq mx+n\}$ (where $y_1$ is the first coordinate of $y \in \mathbb{R}^{n-1}$) with $m$ large, intersects $\partial \Omega$
in a compact hypersurface which is not a graph over a piece of $H=\{y_1=mx+n\}$,
which contradicts Proposition \ref{compact}.\\
Therefore, the intersection of $\Omega$ with any hyperplane normal to the $x$-axis is
nonempty. Define $\Omega_1=\Omega\cap \{x>0\}$
and $\Omega_2=\Omega\cap \{x<0\}$. Note that $\Omega_1$ and
$\Omega_2$ are unbounded. Moreover $\{x=0\} \cap \Omega$ is done of finitely
many open connected domains and the intersection $\{x=0\} \cap \partial \Omega$ is done by more than one point. We can assume that the $x$-axis intersects $\Omega$ transversally. 
If $p$ is the lowest point of $\{x=0\} \cap \partial \Omega$, then we can repeat all the reasoning of the first part of the proof of proposition \ref{halfplane}, with tilted hyperplanes instead of tilted lines, 
getting to the conclusion that there exists a horizontal line $T$ such that $\Omega$ is rotationally symmetric with respect to $T$. Remark that the component $\Sigma(T)$ that can be naturally defined by generalization of the same component in dimension 2, is bounded by the assumption that $\Omega$ is contained in a cylinder. This allows to apply the moving plane argument and completes the proof of the result. \proofend
%
%\begin{figure}[!ht]
%\centering
%{\scalebox{.8}{\input{image7.pstex_t}}}
%\label{KKS}\caption{}
%\end{figure}
%

\section{Concavity properties for double periodic domains}

In this paragraph we deal with 2-dimensional domains that are double periodic, i.e., domains in $\mathbb{R}^{2}$ whose closure is represented by a compact region in the quotient space modulo two linearly independent translations, and where it is possible to solve problem (\ref{pr_bis}). In order to simplify the notation we will consider an open connected domain of a flat torus $T^2 = \mathbb{R}^2\,/\, \langle v_1,v_2\rangle$, where $v_1$ and $v_2$ are two linearly independent vectors. The closure of the connected components of the universal covering of such a domain can be either compact (this case is not interesting because we know that the only bounded domains where it is possible to solve problem (\ref{pr_bis}) are balls), or domains that are periodic in one direction, or a double periodic domain (in this case there exists only one connected component in the covering $\mathbb{R}^{2}$).

We will prove now Theorem \ref{T3}. Before proving the theorem, we want to remark that condition (\ref{hyp}), when $f(u)=\lambda\,u$, becomes
\[
\max_{\overline{\Omega}}\,u <\frac{ |\alpha|}{\sqrt{\lambda}}. 
\]
and for example $\max_{\overline{\Omega}}\,u = |\alpha|/\sqrt{\lambda}$ if $\Omega$ is a the strip $\left(0,\pi/\sqrt{\lambda}\right) \times \mathbb{R}$. %If $f$ is the Allen-Cahn second member, i.e. $f(u)=u-u^3$, and then condition (\ref{hyp}) becomes

The periodicity of $u$ is used just to guarantee that the differential
expressions considered in this section attain their maximum. These
expressions and their connections  with the maximum principle can be found
in the book of R. Sperb \cite{Sperb}, Chapter X.
For constant mean curvature surfaces in the Euclidean space, related
results where proved by A. Ros and H. Rosenberg in \cite{Ros-Ros}. See
\cite{manza} for other ambient spaces.

\medskip

\textit{Proof of Theorem \ref{T3}.} Let us define the operator
\[
P(x) = |\nabla\, u(x)|^2 + 2\, \int_0^{u(x)}f(s)\,\textnormal{d}s
\]
where $P : \Omega \to \mathbb{R}$. Let us denote the coordinates of $\mathbb{R}^2$ by $x=(x_1,x_2)$, and partial derivatives by a comma followed by a subscript, i.e. the partial derivative of a function $u$ with respect to the coordinate $x_i$ will be written as $u_{,i}$ and second partial derivative with respect to the coordinate $x_i$ and $x_j$ will be written as $u_{,ij}$. Moreover we use the standard summation convention. We have
\begin{equation}\label{pi}
P_{,i} = 2\, u_{,ji}\,u_{,j} + 2\, f(u)\, u_{,i}
\end{equation}
and 
\[
\Delta P = P_{,ii} = 2\, u_{,ji}\, u_{,ji} + 2\, u_{,jii}\, u_{,j} + 2\, f'(u)\,|\nabla\, u|^{2} + 2\, f(u)\,\Delta\,u
\]
Using the equation $\Delta\,u +f(u)=0$ and its derivation
\[
u_{,jii} = u_{,iij} = -f'(u)\,u_{,j}
\]
we get
\begin{equation}\label{deltap}
\Delta P = 2\, u_{,ji}\, u_{,ji}  - 2\, f(u)^2
\end{equation}
In order to eliminate the term $2\, u_{,ji}\, u_{,ji}$ we use an identity valid only for two variable functions. In fact, if $v$ is a $C^2$ function of two real variables, an explicit computation shows that
\begin{equation}\label{2var}
|\nabla\, v|^2\, v_{,ij}\,v_{,ij} = |\nabla\, v|^2\, (\Delta\, u)^2 + 2\, v_{,i}\,v_{,ik}\,v_{,j}\,v_{,jk} - 2\, (\Delta\,v)\, v_{,i}\,v_{,j}\,v_{,ij}\,
\end{equation}
Let us define 
\[
L_i = -P_{,i} + 2\,f(u)\,u_{,i}
\]
Using (\ref{pi}), (\ref{deltap}) and (\ref{2var}) we get
\[
\Delta P + \frac{L_i\,P_{,i}}{|\nabla\, u|^2} = 0
\]
The maximum principle can be applied to $P$ at any point $x$ where $\nabla\, u \neq 0$, then, unless $P$ is constant, one of the following two events occurs:
\begin{enumerate}
\item the maximum of $P$ is attaint at $\partial \Omega$, or
\item the maximum of $P$ is attaint at a point $x_0 \in \Omega$ where $\nabla\, u(x_0) = 0$.
\end{enumerate}
We remark that at $\partial \Omega$ we have 
\[
P(x)=\alpha^2
\] 
and, by our hypothesis, at a point $x_0$ where $\nabla\, u(x_0) = 0$ we have
\[
P(x) %\leq 
< \alpha^2
\]
We conclude that $P$ attains its maximum at $\partial \Omega$ and then 
\[
P(x) < \alpha^2
\]
for all $x \in \Omega$ and 
\[
P(x) = \alpha^2
\]
for all $x \in \pp\Omega$.
%\begin{enumerate}
%\item $P(x) = \alpha^2$ for all $x \in \overline \Omega$, or
%\item $P(x) < \alpha^2$ for all $x \in \Omega$ and $P(x) = \alpha^2$ for all $x \in \pp\Omega$.
%\end{enumerate}
So the maximum principle implies the following condition on the normal derivative of $P$:
\begin{equation}\label{pnu}
%\langle \nabla P, \nu \rangle = 0\,\,\textnormal{for all}\,\,x \in \pp \Omega\,\,\textnormal{, or}\,\,
\langle \nabla P, \nu \rangle > 0\,\,\textnormal{for all}\,\,x \in \pp \Omega.
\end{equation} 
Our aim is now to calculate the normal derivative of $P$ in order to make explicit the curvature of $\partial \Omega$. If $t$ is the unit tangent vector about $\pp \Omega$, and we use the same notation as above for derivatives, then at $\pp \Omega$ we have
\begin{eqnarray*}
\langle \nabla P, \nu \rangle & = & 2\, u_{,\nu\nu}\,u_{,\nu} + 2\,f(u)\,u_{,\nu}\\
& = & 2\,u_{,\nu}\, (\Delta u - u_{,tt} + f(u) )\\
& = & - 2\, u_{,\nu}\,u_{,tt}\\
& = & - 2\, \alpha\, u_{,tt}
\end{eqnarray*}
From (\ref{pnu}), and recalling that $\alpha$ is negative, we obtain
\begin{equation}\label{utt}
%u_{,tt} \equiv 0\,\,\textnormal{at}\,\,\pp \Omega\,\,\textnormal{, or}\,\,
u_{,tt} > 0\,\,\textnormal{at}\,\,\pp \Omega.
\end{equation}
Let now $\gamma : S^1 \to \mathbb{R}^2$ be the arclength parametrization of a connected component of $\pp \Omega$ and let $k$ be its curvature with respect to the outward pointing unit normal vector $\nu$. As $u$ is equal to 0 at $\pp \Omega$, we have
\begin{eqnarray*}
0 &=& \langle \nabla u, \gamma' \rangle \\
& = & u_{,tt} + \langle \nabla u, \gamma'' \rangle\\
& = & u_{,tt} + k\, \langle \nabla u, \nu \rangle\\
& = & u_{,tt} + \alpha\,k
\end{eqnarray*}
From (\ref{utt}), and recalling that $\alpha$ is negative, we obtain
\[
%k \equiv 0\,\,\textnormal{at}\,\,\pp \Omega\,\,\textnormal{, or}\,\,
k > 0\,\,\textnormal{at}\,\,\pp \Omega
\]
%In the first case $\Omega$ is a strip, in the second 
i.e., $T^2 \backslash \Omega$ is strictly convex. This completes the proof of the result.

\enddocument